\documentclass[10pt]{amsart}

\usepackage[all]{xy}
\usepackage{latexsym}
\usepackage{amssymb}
\usepackage{a4wide}
\usepackage[utf8]{inputenc}
\usepackage[T1]{fontenc}
\usepackage{dsfont}
\usepackage{pdfsync}

\swapnumbers
\newtheorem{thm}[subsubsection]{Theorem}

\newtheorem{prop}[subsubsection]{Proposition}
\newtheorem{lem}[subsubsection]{Lemma}
\newtheorem{thmAppendix}[subsection]{Theorem}

\theoremstyle{definition}
\newtheorem*{rem}{\sc Remark}

\newtheorem*{pf}{\sc Proof}

\newtheorem*{ex}{\sc Examples}
\newtheorem*{ex1}{\sc Example}

\newcommand{\B}{\mathrm{B}}
\newcommand{\C}{\mathcal{C}}

\newcommand{\F}{\mathcal{F}}

\newcommand{\Po}{\mathcal{P}}
\newcommand{\Poa}{\Po^{\ash}}
\newcommand{\Pos}{\Scc \Po^{\ash}}
\newcommand{\Cs}{\Scc \C}
\newcommand{\So}{\mathbb{S}}
\newcommand{\Sc}{\mathcal{S}}
\newcommand{\Scc}{{\mathcal{S}^{-1}}}
\newcommand{\K}{\mathbb{K}}
\newcommand{\g}{\mathfrak{g}}

\newcommand{\cqfd}{\ \hfill \square}

\newcommand{\As}{\mathcal{A}s}
\newcommand{\Lie}{\mathcal{L}ie}
\newcommand{\Com}{\mathcal{C}om}

\newcommand{\ash}{\textrm{!`}}
\newcommand{\epi}{\twoheadrightarrow}
\newcommand{\mono}{\rightarrowtail}

\newcommand{\qiso}{\xrightarrow{\sim}}

\author{Joan Mill\`es}
\title{The Koszul complex is the cotangent complex}

\begin{document}

\maketitle

\begin{abstract}
We extend the Koszul duality theory of associative algebras to algebras over an operad. Recall that in the classical case, this Koszul duality theory relies on an important chain complex: the Koszul complex. We show that the cotangent complex, involved in the cohomology theory of algebras over an operad, generalizes the Koszul complex.
\end{abstract}

\section*{Introduction}

The Koszul duality theory of associative algebras first appeared in the work of \cite{Priddy}. It has then been extended to operads \cite{GetzlerJones, GinzburgKapranov} and to properads \cite{Vallette2}. The Koszul duality theory has a wide range of applications in mathematics: BGG correspondence \cite{BGG}, equivariant cohomology \cite{GKMP}, and homotopy algebras \cite{GinzburgKapranov, GetzlerJones}. For even more applications, we refer to the introduction of \cite{PolishchukPositselski}. Following the ideas of \cite{Quillen}, one defines a cohomology theory associated to any type of algebras \cite{Hinich, GoerssHopkins} and the Koszul duality theory of operads provides explicit chain complexes which allow us to compute it \cite{Milles}.\\

The Koszul duality theory of associative algebras is based on a chain complex, called the \emph{Koszul complex}, built on the tensor product of chain complexes. The notions of algebras, operads and properads are ``associative'' notions in the sense that they are all monoids in a monoidal category. This makes the generalization of the \emph{Koszul complex} to operads and properads possible. To define a Koszul duality theory for algebras over an operad $\Po$, we have to find a good generalization for this Koszul complex in a non-associative setting.\\

In \cite{Milles}, we studied the André-Quillen cohomology theory of algebra over an operad. The latter is represented by a chain complex, called the \emph{cotangent complex} of a $\Po$-algebra $A$. Thanks to the Koszul duality theory of operads, we made a representation of the cotangent complex explicit. In this paper, we prove that the cotangent complex gives a good generalization of the Koszul complex in the sense that we get an algebraic twisting morphisms fundamental theorem (Theorem \ref{atmft}) and a Koszul criterion (Theorem \ref{algKoszulcriterion}).\\

When the $\Po$-algebra $A$ is \emph{quadratic} or \emph{monogene}, we introduce a Koszul dual coalgebra $A^{\ash}$. The Koszul criterion provides a way to test whether the Koszul dual coalgebra $A^{\ash}$ is a good space of \emph{syzygies} to resolve the $\Po$-algebra $A$. If applicable, the Koszul complex, which is thus a representation of the cotangent complex, is a ``small'' chain complex allowing to compute the cohomology theory of the $\Po$-algebra $A$.\\

Retrospectively, the present Koszul duality theory applied to associative algebras gives the Koszul duality theory of associative algebras originally defined by Priddy \cite{Priddy}. For commutative algebras, resp. Lie algebras, the present Koszul duality theory provides Sullivan minimal models, resp. Quillen models. An example is given by the commutative algebra of the cohomology groups of the complement of an hyperplane arrangement. It is given by the Orlik-Solomon algebra and, in the quadratic case, the Koszul dual algebra is the holonomy Lie algebra defined by Kohno \cite{Yuzvinsky, PapadimaYuzvinsky, Kohno, Kohno2}. More generally, this theory applies to give all the rational homotopy groups of formal spaces whose cohomology groups forms a Koszul (quadratic) algebra. Moreover, an associative algebra $A$, eventually commutative, is an example of operad and an ``algebra'' over this operad is a $A$-module. The present Koszul duality theory in this case gives the Koszul duality theory of modules which provides a good candidate for the syzygies of an $A$-module \cite{PolishchukPositselski, Eisenbud}.\\

The following table gives a summary.
$$\begin{tabular}{l||p{3.1cm}|p{3.4cm}|}
\cline{2-3}
& \multicolumn{2}{c|}{\bf Koszul duality theory}\\
\hline \hline
\multicolumn{1}{|p{3.5cm}||}{{\bf Monoids}\hspace{2cm} Section 1} & Associative algebras \cite{Priddy} & Operads \cite{GinzburgKapranov, GetzlerJones}\\
\hline
\multicolumn{1}{|p{3.5cm}||}{{\bf Representations}\hspace{2cm} Sections 2, 3 and 4} & Modules \hspace{1cm} \cite{PolishchukPositselski} & $\Po$-algebras\hspace{1.5cm} {\it Goal} of the paper\\
\hline
\end{tabular}$$\\

We recall in the first section the Koszul duality theory for associative algebras and operads. We recall the results of \cite{GetzlerJones} on twisting morphisms for $\Po$-algebras and we prove the algebraic twisting morphisms fundamental theorem in Section 2. We extend the results of \cite{GinzburgKapranov} on quadratic $\Po$-algebras and state the Koszul criterion in Section 3. The section 4 is devoted to the applications of the Koszul duality theory for $\Po$-algebras. We prove a comparison Lemma for the twisted tensor product in the framework of $\Po$-algebras in Appendix A.\\

In all this paper, $\K$ is a field of characteristic $0$, however most constructions work over a ring. Moreover, we assume that all the chain complexes are non-negatively graded.

\tableofcontents

\section{Twisting morphisms for associative algebras and operads}\label{twbar}

The notion of \emph{twisting morphism} (or twisting cochain) for associative algebras, introduced by Cartan \cite{Cartan} and Brown \cite{Brown}, was generalized to operads by Getzler and Jones \cite{GetzlerJones}. We recall their definitions and the fact that the induced bifunctor is represented by the bar and the cobar constructions. We also recall the definition of the \emph{twisted tensor product} and the \emph{twisted composition product} and the notion of \emph{(operadic) Koszul morphism}. We gives the twisting morphisms fundamental theorems and the Koszul criteria. We refer to the book of Loday and Vallette \cite{LodayVallette} for a complete exposition.

\subsection{Twisting morphisms for associative algebras}\label{casas}

Let $(A,\, \gamma_{A},\, d_{A})$ be a \emph{differential graded associative algebra}, \emph{dga algebra} for short, and let $(C,\, \Delta_{C},\, d_{C})$ be a \emph{differential graded coassociative coalgebra}, \emph{dga coalgebra} for short. We associate to $C$ and $A$ the \emph{dg convolution algebra} $(\mathrm{Hom}(C,\, A),\, \star,\, \partial)$ where $\star$ and $\partial$ are defined as follows
$$\left\{\begin{array}{l}
f \star g : C \xrightarrow{\Delta_{C}} C \otimes C \xrightarrow{f \otimes g} A \otimes A \xrightarrow{\gamma_{A}} A,\\
\partial(f) := d_{A} \circ f - (-1)^{|f|}f \circ d_{C},
\end{array} \right.$$
for $f,\, g \in \mathrm{Hom}(C,\, A)$. The associative product $\star$ on $\mathrm{Hom}(C,\, A)$ induces a Lie bracket and the solutions of degree $-1$ of the \emph{Maurer-Cartan equation}
$$\partial (\alpha) + \alpha \star \alpha = \partial(\alpha) + \frac{1}{2} [\alpha,\, \alpha] = 0$$
are called \emph{twiting morphisms}. The set of twisting morphisms is denoted by $\mathrm{Tw}(C,\, A)$.

When $A$ is augmented, that is $A \cong \K \oplus \overline{A}$ as dga algebras, we recall the \emph{bar construction on $A$} defined by $\B A := (T^{c}(s\overline{A}),\, d := d_{1} + d_{2})$, where $s$ is the homological suspension, $d_{1}$ is induced by the differential $d_{A}$ on $A$ and $d_{2}$ is the unique coderivation which extends, up to suspension, the restriction of the product $\gamma_{A}$ on $\overline{A}$. When $C$ is coaugmented, that is $C \cong \K \oplus \overline{C}$ as coassociative coalgebras,  we recall dually the \emph{cobar construction on $C$} defined by $\Omega C := (T(s^{-1}\overline{C}),\, d := d_{1} - d_{2})$, where $s^{-1}$ is the homological desuspension, $d_{1}$ is induced by $d_{C}$ and $d_{2}$ is the unique derivation which extends, up to desuspension, the restriction of $\Delta_{C}$ on $\overline{C}$.

When $A$ is augmented and $C$ is coaugmented, we require that the composition of a twisting morphism with the augmentation map, respectively the coaugmentation map, vanishes. These constructions satisfy the following adjunction
$$\begin{array}{c}
\mathrm{Hom_{\mathsf{aug.\ dga\ alg.}}}(\Omega C,\, A) \cong \mathrm{Tw}(C,\, A) \cong \mathrm{Hom_{\mathsf{coaug.\ dga\ coalg.}}}(C,\, \B A)\\
\xymatrix@C=40pt@M=8pt{f_{\alpha} & \ar@{<->}[l] \alpha \ar@{<->}[r] &  g_{\alpha},}
\end{array}$$
when the coalgebra $C$ is \emph{conilpotent} (see \cite{LodayVallette} for a definition).

To a twisting morphism $\alpha$ between a coaugmented coalgebra $C$ and an augmented algebra $A$, we associate the \emph{left twisted tensor product}
$$A \otimes_{\alpha} C := (A \otimes C,\, d_{\alpha} := d_{A\otimes C} - d_{\alpha}^{l}),$$
where $d_{\alpha}^{l}$ is defined by
$$A \otimes C \xrightarrow{id_{A} \otimes \Delta_{C}} A \otimes C \otimes C \xrightarrow{id_{A} \otimes \alpha \otimes id_{C}} A \otimes A \otimes C \xrightarrow{\gamma_{A} \otimes id_{C}} A \otimes C.$$
It is a chain complex since $\alpha$ is a twisting morphism. We refer to \cite{LodayVallette} for the symmetric definition of the \emph{right twisted tensor product}
$$C \otimes_{\alpha} A := (C \otimes A,\, d_{\alpha} := d_{C \otimes A} + d_{\alpha}^{r})$$
and one gets the \emph{twisted tensor product}
$$A \otimes_{\alpha} C \otimes_{\alpha} A := (A \otimes C \otimes A,\, d_{\alpha} := d_{A \otimes C \otimes A} - d_{\alpha}^{l} \otimes id_{A} + id_{A} \otimes d_{\alpha}^{r}).$$
We say that $\alpha$ is a \emph{Koszul morphism} when $A \otimes_{\alpha} C \otimes_{\alpha} A \qiso A$. We denote by $\mathrm{Kos}(C,\, A)$ the set of Koszul morphisms.

\begin{ex}
To $id_{\B A}$ and $id_{\Omega C}$ correspond two universal twisting morphisms $\pi : \B A \rightarrow A$ and $\iota : C \rightarrow \Omega C$. They are examples of Koszul morphisms, that is $A \otimes_{\pi} \B A \otimes_{\pi} A \qiso A$ and $\Omega C \otimes_{\iota} C \otimes_{\iota} \Omega C \qiso \Omega C$.
\end{ex}

Later, we will need an extra grading, called \emph{weight grading}, which differs from the homological grading. We refer to the first chapter of \cite{LodayVallette} for more details about this and for a definition of connected wdga algebras and connected wdga coalgebras. This adjunction satisfies the following property.
\begin{thm}[Twisting morphisms fundamental theorem]\label{tmft}
Let $A$ be a connected wdga algebra and let $C$ be a connected wdga coalgebra. For any twisting morphism $\alpha : C \rightarrow A$, the following assertions are equivalent:
\begin{enumerate}
\item The twisting morphism $\alpha$ is Koszul, that is $A \otimes_{\alpha} C \otimes_{\alpha} A \qiso A$;
\item The left twisted tensor product is acyclic, that is $A \otimes_{\alpha} C \qiso \K$;
\item The right twisted tensor product is acyclic, that is $C \otimes A \qiso \K$;
\item The morphism of dga algebras $f_{\alpha} : \Omega C \rightarrow A$ is a quasi-isomorphism;
\item The morphism of dga coalgebras $g_{\alpha} : C \rightarrow \B A$ is a quasi-isomorphism.
\end{enumerate}
\end{thm}

\begin{pf}
A proof of the equivalences $(2) \Leftrightarrow (3) \Leftrightarrow (4) \Leftrightarrow (5)$ can be found in \cite{LodayVallette} and comes from \cite{Brown}. We show the equivalence $(1) \Leftrightarrow (2)$ in the more general case of operads, see Theorem \ref{otmft}.
$\cqfd$
\end{pf}

Let $(V,\, S)$ be a \emph{quadratic data}, that is a graded vector space $V$ and a graded subspace $R \subseteq V \otimes V$. A \emph{quadratic algebra} is an associative algebra $A(V,\, S)$ of the form $T(V)/(S)$. Dually the \emph{quadratic coalgebra} $C(V,\, S)$ is by definition the sub-coalgebra of the cofree coalgebra $T^{c}(V)$ which is \emph{universal} among the sub-coalgebras $C$ of $T^{c}(V)$ such that the composite
$$C \mono \F^{c}(V) \epi \F^{c}(V)^{(2)}/S$$
is $0$. The word ``universal'' means that for any such coalgebra $C$, there exists a unique morphism of coalgebras $C \rightarrow C(V,\, S)$ such that the following diagram commutes
$$\xymatrix@M=6pt@R=12pt{C \ar@{>->}[dr] \ar[r] & C(V,\, S) \ar@{>->}[d] \\
& T^{c}(V).}$$

To a quadratic data $(V,\, S)$, we associate the \emph{Koszul dual coalgebra of $A$} given by $A^{\ash} := C(sV,\, s^{2}S)$ where $s$ is the homological suspension. We associate to the coalgebra $C(sV,\, s^{2}S)$ and to the algebra $A(V,\, S)$ the twisting morphism $\varkappa$ defined by
$$\varkappa : A^{\ash} = C(sV,\, s^{2}S) \epi sV \cong V \mono A(V,\, S) = A.$$
The equality $\partial(\varkappa) + \varkappa \star \varkappa = 0$ follows from the coassociativity of $\Delta_{A^{\ash}}$, the associativity of $\gamma_{A}$ and to the fact that ${A^{\ash}}^{(2)} = s^{2}S$ and $A^{(2)} = V^{\otimes 2}/S$.

The weight grading comes from the graduation in $T(V)$ and $T^{c}(V)$ by the number of generators in $V$. The universality of the twisting morphisms $\pi$ and $\iota$ provides an inclusion of coalgebras $g_{\varkappa} : A^{\ash} \mono \B A$ and a surjection of algebras $f_{\varkappa} : \Omega A^{\ash} \epi A$. The twisting morphisms fundamental theorem writes:
\begin{thm}[Koszul criterion, \cite{Priddy}]\label{Koszulcriterion}
Let $(V,\, S)$ be a quadratic data. Let $A := A(V,\, S)$ be the associated quadratic algebra and let $A^{\ash} := C(sV,\, s^{2}S)$ its Koszul dual coalgebra. The following assertions are equivalent:
\begin{enumerate}
\item The twisting morphism $\varkappa$ is Koszul, that is $A \otimes_{\varkappa} A^{\ash} \otimes_{\varkappa} A \qiso A$;
\item The Koszul complex $A \otimes_{\varkappa} A^{\ash}$ is acyclic, that is $A \otimes_{\varkappa} A^{\ash} \qiso \K$;
\item The Koszul complex $A^{\ash} \otimes_{\varkappa} A$ is acyclic, that is $A^{\ash} \otimes_{\varkappa} A \qiso \K$;
\item The morphism of dga algebras $f_{\varkappa} : \Omega A^{\ash} \rightarrow A$ is a quasi-isomorphism;
\item The morphism of dga coalgebras $g_{\varkappa} : A^{\ash} \rightarrow \B A$ is a quasi-isomorphism.
\end{enumerate}
\end{thm}

Priddy \cite{Priddy} called the chain complexes $A^{\ash} \otimes_{\varkappa} A$, resp. $A \circ_{\varkappa} A^{\ash}$, the \emph{Koszul complexes}. We extend this definition and we call also $A \otimes_{\varkappa} A^{\ash} \otimes_{\varkappa} A$ the Koszul complex. An algebra is called a \emph{Koszul algebra} when the twisting morphism $\varkappa : A^{\ash} \rightarrow A$ is a Koszul morphism, that is $A \circ_{\varkappa} A^{\ash} \otimes_{\varkappa} A \qiso A$. The previous Koszul criterion 
shows that this definition is equivalent to the classical one.

\subsection{$\So$-module and operad}

We recall the definition of an \emph{$\So$-module} and of an \emph{operad}. For a complete exposition of the concepts of this section, we refer to the books \cite{LodayVallette} and \cite{MarklShniderStasheff}.\\

An \emph{$\So$-module} $M$ is a collection of dg modules $\{ M(n)\}_{n \geq 0}$ endowed with an action of the group $\So_{n}$ of permutations on $n$ elements. Let $M,\, M',\, N$ and $N'$ be $\So$-modules. We recall the definition of the composition product $\circ$
$$(M\circ N) (n) := \bigoplus_{k \geq 0} M(k) \otimes_{\So_k} \left(\bigoplus_{i_1 + \cdots + i_k = n}
\mathrm{Ind}_{\So_{i_{1}}\times \cdots \times \So_{i_{k}}}^{\So_{n}} \big(N(i_1) \otimes \cdots \otimes N(i_k)\big)\right).$$
The unit for the monoidal product is $I := (0,\, \K,\, 0,\, \ldots)$. Notice that $\circ$ is not linear on the right hand side. As a consequence, we define the right linear analog $M \circ (N;\, N')$ of the composition product, linear in $M$ and in $N'$, by the following formula
\begin{eqnarray*}
&& M \circ (N;\, N')(n) :=\\
&& \bigoplus_{k \geq 0} M(k) \otimes_{\So_k} \left(\bigoplus_{i_1 + \cdots + i_k = n} \bigoplus_{j = 1}^{k} \mathrm{Ind}_{\So_{i_{1}}\times \cdots \times \So_{i_{k}}}^{\So_{n}} (N(i_1) \otimes \cdots \otimes \underbrace{N'(i_{j})}_{j^{th} \textrm{ position}} \otimes \cdots \otimes N(i_k))\right).
\end{eqnarray*}
Let $f : M \rightarrow M'$ and $g : N \rightarrow N'$ be maps of $\So$-modules. We denote by $\circ'$ the \emph{infinitesimal composite of morphisms}:
$$f \circ' g : M \circ N \rightarrow M' \circ (N,\, N')$$
defined by
$$\sum_{j=1}^k f \otimes (id_{N} \otimes \cdots \otimes \underbrace{g}_{j^{th} \textrm{ position}} \otimes \cdots \otimes id_{N}).$$
The differential on $M\circ N$ is given by $d_{M\circ N} := d_M \circ id_N + id_M \circ' d_N$. The term $id_{M} \circ' d_{N}$ goes normally to $M \circ (N;\, N)$ but we assume it composed with the projection $M \circ (N;\, N) \epi M \circ N$. We define the \emph{infinitesimal composite product} $M\circ_{(1)} N$ by $M\circ (I;\, N)$, which is linear in $M$ and in $N$. Moreover we denote by $f\circ_{(1)} g$ the map $f \circ (id_{I},\, g) : M\circ_{(1)} N \rightarrow M'\circ_{(1)} N'$.

The category of dg $\So$-modules $(\So \textrm{-mod},\, \circ,\, I)$ is a monoidal category. A monoid $(\Po,\, \gamma,\, u)$ in this category is called a \emph{differential graded operad}, \emph{dg operad} for short. Dually, a comonoid $(\C,\, \Delta,\, \eta)$ in this category is called a \emph{dg cooperad}.

Let $(\Po,\, \gamma,\, u)$ be a dg operad and $(\C,\, \Delta,\, \eta)$ be a dg cooperad. We denote by $\gamma_{(1)}$ the \emph{partial product} of the operad $\Po$
$$\Po \circ_{(1)} \Po \mono \Po \circ \Po \xrightarrow{\gamma} \Po.$$
Dually, we denote by $\Delta_{(1)}$ the \emph{partial coproduct} of the cooperad $\C$
$$\C \xrightarrow{\Delta} \C \circ \C \epi \C \circ_{(1)} \C.$$

The free operad $\F(E)$ on a $\So$-module $E$ is given by all the trees whose vertices of arity $n$ are indexed by elements in $E(n)$. The product is given by the grafting of trees. Dually, the free cooperad $\F^{c}(E)$ has the same underlying $\So$-module as $\F(E)$ and the coproduct is given by the splitting of trees.

\subsection{Operadic twisting morphism}

We recall from \cite{GetzlerJones, GinzburgKapranov} the definitions of \emph{operadic twisting morphism}, \emph{left twisted composition product} and \emph{bar and cobar construction} in the setting of operads. We state the operadic twisting morphisms fundamental theorem and recall the notion of \emph{Koszul operad}.\\

Let $\Po$ be a dg operad and $\C$ be a dg cooperad. We recall from Chapter 6 of \cite{LodayVallette} the \emph{dg convolution PreLie algebra} $(\mathrm{Hom_{\So\textrm{-mod}}}(\C,\, \Po),\, \star,\, \partial)$ where
$$\left\{ \begin{array}{l}
f\star g : \C \xrightarrow{\Delta_{(1)}} \C \circ_{(1)} \C \xrightarrow{f\circ_{(1)}g} \Po \circ_{(1)} \Po \xrightarrow{\gamma_{(1)}} \Po,\\
\partial(f) := d_{\Po} \circ f - (-1)^{|f|} f \circ d_{\C},
\end{array} \right.$$
for $f,\, g \in \mathrm{Hom}(\C,\, \Po)$. The PreLie product $\star$ induces a Lie bracket on $\mathrm{Hom}(\C,\, \Po)$ and an \emph{operadic twisting morphism} is a map $\alpha : \C \rightarrow \Po$ of degree $-1$, solution of the \emph{Maurer-Cartan equation}
$$\partial (\alpha) + \alpha \star \alpha = \partial(\alpha) + \frac{1}{2} [\alpha,\, \alpha] = 0.$$
We denote the set of operadic twisting morphisms from $\C$ to $\Po$ by Tw$(\C,\, \Po)$.

In \cite{GetzlerJones, GinzburgKapranov}, the authors extend the bar construction and the cobar construction of associative algebras/coalgebras to operads/cooperads. When $\Po$ is an \emph{augmented operad}, that is $\Po \cong I \oplus \overline{\Po}$, we have $\B \Po := (\F^{c}(s\overline{\Po}),\, d := d_{1} + d_{2})$ where $\F^{c}(s\overline{\Po})$ is the free cooperad on the homological suspension of $\overline{\Po}$, $d_{1}$ is the unique coderivation which extends, up to suspension, the differental $d_{\overline{\Po}}$ and $d_{2}$ is the unique coderivation which extends, up to suspension, the restriction of the partial product $\overline{\gamma}_{(1)} : \overline{\Po} \circ_{(1)} \overline{\Po} \rightarrow \overline{\Po}$. When $\C$ is a \emph{coaugmented cooperad}, that is $\C \cong I \oplus \overline{\C}$ as cooperads, we have $\Omega \C := (\F(s^{-1}\overline{\C}),\, d := d_{1} - d_{2})$, where $\F(s^{-1}\overline{\C})$ is the free operad on the homological desuspension of $\overline{\C}$, $d_{1}$ is the unique derivation which extends, up to desuspension, the differential $d_{\overline{\C}}$ and $d_{2}$ is the unique derivation which extends, up to desuspension, the partial coproduct $\overline{\Delta}_{(1)} : \overline{\C} \rightarrow \overline{\C} \circ_{(1)} \overline{\C}$. When $\Po$ is an augmented operad and $\C$ is a coaugmented cooperad, we require that the composition of an operadic twisting morphism with the augmentation map, respectively the coaugmentation map, vanishes. As for associative algebras, these constructions satisfy the following bar-cobar adjunction.

\begin{thm}[Theorem $2.17$ of \cite{GetzlerJones}]
The functors $\Omega$ and $\B$ form a pair of adjoint functors between the category of conilpotent coaugmented dg cooperads and augmented dg operads. The natural bijections are given by the set of operadic twisting morphisms:
$$\begin{array}{c}
\mathrm{Hom}_{\mathsf{dg\ op.}}(\Omega \C,\, \Po) \cong \mathrm{Tw(\C, \Po)} \cong \mathrm{Hom}_{\mathsf{dg\ coop.}}(\C,\, \B \Po)\\
\xymatrix@C=40pt@M=8pt{f_{\alpha} & \ar@{<->}[l] \alpha \ar@{<->}[r] &  g_{\alpha}.}
\end{array}$$
\end{thm}

\begin{ex}[of operadic twisting morphisms]
\begin{itemize}
\item[]
\item When $\C = \B \Po$, the previous theorem gives a natural operadic twisting morphism $\pi : \B \Po \rightarrow \Po$, associated to $id_{\B \Po}$, which is equal to $\B \Po = \mathcal{F}^{c}(s\overline{\Po}) \epi s \overline{\Po} \xrightarrow{s^{-1}} \overline{\Po} \mono \Po$. This morphism is universal in the sense that each operadic twisting morphism $\alpha : \C \rightarrow \Po$ factorizes through $\pi$
$$\xymatrix{\C \ar[rr]^{\alpha} \ar@{-->}[dr]_{f_{\alpha}} && \Po\\
& \B \Po, \ar[ur]_{\pi}}$$
where $f_{\alpha}$ is a morphism of dg cooperads.
\item When $\Po = \Omega \C$, the previous theorem gives a natural operadic twisting morphism $\iota : \C \rightarrow \Omega \C$, associated to $id_{\Omega \C}$, which is equal to $\C \epi \overline{\C} \xrightarrow{s^{-1}} s^{-1}\overline{\C} \mono \Omega \C = \F (s^{-1}\overline{\C})$. This morphism is universal in the sense that each operadic morphism $\alpha : \C \rightarrow \Po$ factorizes through $\iota$
$$\xymatrix{& \Omega \C \ar@{-->}[dr]^{g_{\alpha}} &\\
\C \ar[ur]^{\iota} \ar[rr]^{\alpha} && \Po,}$$
where $g_{\alpha}$ is a morphism of dg operads.
\end{itemize}
\end{ex}

To an operadic twisting morphism $\alpha$ between a coaugmented cooperad $\C$ and an augmented operad $\Po$, we associate the \emph{left twisted composition product} \cite{Vallette2, LodayVallette}
$$\Po \circ_{\alpha} \C := (\Po \circ \C,\, d_{\alpha} := d_{\Po \circ \C} - d_{\alpha}^{l}),$$
where $d_{\alpha}^{l}$ is given by the composite
$$\Po \circ \C \xrightarrow{id_{\Po} \circ' \Delta} \Po \circ (\C;\, \C \circ \C) \xrightarrow{id_{\Po} \circ (id_{\C},\, \alpha \circ id_{\C})} \Po \circ (\C;\, \Po \circ \C) \cong (\Po \circ_{(1)} \Po) \circ \C \xrightarrow{\gamma_{(1)} \circ id_{\C}} \Po \circ \C.$$
Since $\alpha$ is an operadic twisting morphism, $d_{\alpha}$ is a differential (see \cite{LodayVallette}). We refer to \cite{LodayVallette} for a definition of the \emph{right twisted composite product}
$$\C \circ_{\alpha} \Po := (\C \circ \Po,\, d_{\alpha} := d_{\C \circ \Po} + d_{\alpha}^{r})$$
and one gets the \emph{twisted composite product}
$$\Po \circ_{\alpha} \C \circ_{\alpha} \Po := (\Po \circ \C \circ \Po,\, d_{\alpha} := d_{\Po \circ \C \circ \Po} - d_{\alpha}^{l} \circ id_{\Po} + id_{\Po} \circ d_{\alpha}^{r}).$$
We say that $\alpha$ is an \emph{operadic Koszul morphism} when $\Po \circ_{\alpha} \C \circ_{\alpha} \Po \qiso \Po$. We denote by $\mathrm{Kos}(\C,\, \Po)$ the set of operadic Koszul morphisms.

\begin{lem}[\cite{GetzlerJones, Fresse, Vallette2}]\label{augbarres}
The twisting morphisms $\pi : \B \Po \rightarrow \Po$ and $\iota : \C \rightarrow \Omega \C$ are operadic Koszul morphisms, that is
$$\Po \circ_{\pi} \B \Po \circ_{\pi} \Po \qiso \Po \textrm{ and } \Omega \C \circ_{\iota} \C \circ_{\iota} \Omega \C \qiso \Omega \C.$$
\end{lem}

Sometimes, we need an extra grading, called \emph{weight grading}, which differs from the homological degree. We say that a \emph{weight graded dg $\So$-module}, \emph{wdg $\So$-module} for short, $M$ is \emph{connected} when $M^{(0)} = I$ and $M = I \oplus M^{(1)} \oplus \cdots \oplus M^{(\omega)} \oplus \cdots$. These definitions hold for operads and cooperads. For example, the weight grading on the free operad $\F(E)$ or on the free cooperad $\F^{c}(E)$ is given by the number $\omega$ of vertices and denoted by $\F (E)^{(\omega)}$ or $\F^{c}(E)^{(\omega)}$. This induces a weight grading on each \emph{quadratic operad} $\Po = \F(E)/(R)$, where $R \subset \F(E)^{(2)}$ and on each sub-cooperad of a free cooperad. In the weight graded setting, we assume that the maps preserve the weight grading. For example, a twisting morphism $\alpha$ preserves the weight grading and when the underlying modules are connected, we have $\alpha = \alpha^{(\geq 1)}$.

As for associative algebras, the twisted composite products and the bar and cobar constructions satisfy the following property.
\begin{thm}[Operadic twisting morphisms fundamental theorem, \cite{LodayVallette}]\label{otmft}
Let $\Po$ be a connected wdg operad and let $\C$ be a connected wdg cooperad. For any operadic twisting morphism $\alpha : \C \rightarrow \Po$, the following assertions are equivalent:
\begin{enumerate}
\item The twisting morphism $\alpha$ is Koszul, that is $\Po \circ_{\alpha} \C \circ_{\alpha} \Po \qiso \Po$;
\item The left twisted composite product $\Po \circ_{\alpha} \C$ is acyclic, that is $\Po \circ_{\alpha} \C \qiso I$;
\item The right twisted composite product $\C \circ_{\alpha} \Po$ is acyclic, that is $\C \circ_{\alpha} \Po \qiso I$;
\item The morphism of dg operads $f_{\alpha} : \Omega \C \rightarrow \Po$ is a quasi-isomorphism;
\item The morphism of dg cooperads $g_{\alpha} : \C \rightarrow \B \Po$ is a quasi-isomorphism.
\end{enumerate}
\end{thm}

\begin{pf}
A proof of the equivalences $(2) \Leftrightarrow (3) \Leftrightarrow (4) \Leftrightarrow (5)$ can be find in \cite{LodayVallette}. We show the equivalence $(1) \Leftrightarrow (2)$.

$(1) \Rightarrow (2)$: The Koszul complex $\Po \circ_{\alpha} \C$ is equal to the relative composite product $(\Po \circ_{\alpha} \C \circ_{\alpha} \Po) \circ_{\Po} I$ which is defined by the short exact sequence
$$0 \rightarrow (\Po \circ_{\alpha} \C \circ_{\alpha} \Po) \circ \Po \circ I \rightarrow (\Po \circ_{\alpha} \C \circ_{\alpha} \Po) \circ I \rightarrow (\Po \circ_{\alpha} \C \circ_{\alpha} \Po) \circ_{\Po} I \rightarrow 0.$$
Since we work over a field of characteristic zero, the ring $\K [\So_{n}]$ is semi-simple by Maschke's theorem, that is every $\K [\So_{n}]$-module is projective. So the K\"unneth formula implies $\mathrm{H}_{\bullet}((\Po \circ_{\alpha} \C \circ_{\alpha} \Po) \circ \Po \circ I) \cong \mathrm{H}_{\bullet}(\Po) \circ \mathrm{H}_{\bullet}(\Po) \circ I$. Moreover $\mathrm{H}_{\bullet}((\Po \circ_{\alpha} \C \circ_{\alpha} \Po) \circ I) \cong \mathrm{H}_{\bullet}(\Po) \circ I$. Finally, this gives $\mathrm{H}_{\bullet}(\Po \circ_{\alpha} \C) \cong \mathrm{H}_{\bullet}((\Po \circ_{\alpha} \C \circ_{\alpha} \Po) \circ_{\Po} I) \cong \mathrm{H}_{\bullet}(\Po) \circ_{\mathrm{H}_{\bullet}(\Po)} I \cong I$.

$(2) \Rightarrow (1)$: We define a filtration $F_{p}$ on $\Po \circ_{\alpha} \C \circ_{\alpha} \Po$ by
$$F_{p}(\Po \circ_{\alpha} \C \circ_{\alpha} \Po) := \oplus_{\omega \leq p}(\Po \circ_{\alpha} \C)^{(\omega)} \circ_{\alpha} \Po.$$
The differential on $\Po \circ_{\alpha} \C \circ_{\alpha} \Po$ is given by $d_{\alpha} := d_{\Po \circ \C \circ \Po} - d_{\alpha}^{l} \circ id_{\Po} + id_{\Po} \circ d_{\alpha}^{r}$ and satisfies
$$\left\{ \begin{array}{ccl}
d_{\Po \circ \C \circ \Po} & : & F_{p} \rightarrow F_{p}\\
d_{\alpha}^{l} & : & F_{p} \rightarrow F_{p}\\
d_{\alpha}^{r} & : & F_{p} \rightarrow F_{p-1}.
\end{array} \right.$$
So the filtration is a filtration of chain complexes. Moreover, the filtration is exhaustive and bounded below. We can apply the classical theorem of convergence of spectral sequences (Theorem 5.5.1 of \cite{Weibel}) and we get that the induced spectral sequence $E_{p,\, q}^{\bullet}$ converges to the homology of $\Po \circ_{\alpha} \C \circ_{\alpha} \Po$. We consider the trivial chain complex filtration on $\Po$, that is $F_{p}\Po := \Po$ for all $\Po$, so that the map $\Po \circ_{\alpha} \C \circ_{\alpha} \Po \rightarrow \Po$ respects the filtration. This last map induces a map on the $E^{1}$-pages which is an isomorphism since $E_{p,\, q}^{1} \cong \mathrm{H}_{\bullet}(\Po \circ_{\alpha} \C) \circ \mathrm{H}_{\bullet}(\Po) \cong \mathrm{H}_{\bullet}(\Po) \cong E_{0,\, q}^{1}$ by the K\"unneth formula. The convergence of the spectral sequences concludes the proof.
$\cqfd$
\end{pf}

\subsubsection{\bf Koszul operad}\label{Koszuloperad}

An \emph{operadic quadratic data} $(E,\, R)$ is a graded $\So$-module $E$ and a graded sub-$\So$-module $R \subset \F(E)^{(2)}$. The quotient $\Po(E,\, R) := \F(E)/(R)$ is called a \emph{quadratic operad}. Dually the \emph{quadratic cooperad} $\C(E,\, R)$ is by definition the sub-cooperad of the cofree cooperad $\F^{c}(E)$ which is \emph{universal} among the sub-cooperads $\C$ of $\F^{c}(E)$ such that the composite
$$\C \mono \F^{c}(E) \epi \F^{c}(E)^{(2)}/R$$
is $0$. The word ``universal'' means that for any such cooperad $\C$, there exists a unique morphism of cooperads $\C \rightarrow \C(E,\, R)$ such that the following diagram commutes
$$\xymatrix@M=6pt@R=12pt{\C \ar@{>->}[dr] \ar[r] & \C(E,\, R) \ar@{>->}[d] \\
& \F^{c}(E).}$$

To a quadratic data $(E,\, R)$, we associate the \emph{Koszul dual cooperad of $\Po$} given by $\Poa := \C(sE,\, s^{2}R)$ where $s$ is the homological suspension.

\begin{ex1}[of operadic twisting morphism]
When $\C = \Po^{\ash}$, the map $\kappa : \Po^{\ash} \mono \B \Po \xrightarrow{\pi} \Po$ is an operadic twisting morphism. It is equal to $\Po^{\ash} \epi sE \xrightarrow{s^{-1}} E \mono \Po$.
\end{ex1}

The universality of the twisting morphisms $\pi$ and $\iota$ provides an inclusion of cooperads $g_{\kappa} : \Poa \mono \B \Po$ and a surjection of operads $f_{\kappa} : \Omega \Poa \epi \Po$. The operadic twisting morphisms fundamental theorem writes:
\begin{thm}[Koszul criterion, \cite{LodayVallette}]\label{opKoszulcriterion}
Let $(E,\, R)$ be an operadic quadratic data. Let $\Po := \Po(E,\, R)$ be the associated quadratic operad and let $\Poa := \C(sE,\, s^{2}R)$ its Koszul dual cooperad. The following assertions are equivalent:
\begin{enumerate}
\item The operadic twisting morphism $\kappa$ is Koszul, that is $\Po \circ_{\kappa} \Poa \circ_{\kappa} \Po \qiso \Po$;
\item The Koszul complex $\Po \circ_{\kappa} \Poa$ is acyclic, that is $\Po \circ_{\kappa} \Poa \qiso I$;
\item The Koszul complex $\Poa \circ_{\kappa} \Po$ is acyclic, that is $\Poa \circ_{\kappa} \Po \qiso I$;
\item The morphism of dg operads $f_{\kappa} : \Omega \Poa \rightarrow \Po$ is a quasi-isomorphism;
\item The morphism of dg cooperads $g_{\kappa} : \Poa \rightarrow \B \Po$ is a quasi-isomorphism.
\end{enumerate}
\end{thm}

The chain complexes $\Po \circ_{\kappa} \Poa \circ_{\kappa} \Po$, resp. $\Poa \circ_{\kappa} \Po$, resp. $\Po \circ_{\kappa} \Poa$, are called the \emph{Koszul complexes}. An operad is called a \emph{Koszul operad} when the operadic twisting morphism $\kappa : \Poa \rightarrow \Po$ is an operadic Koszul morphism, that is $\Po \circ_{\kappa} \Poa \circ_{\kappa} \Po \qiso \Po$.

\section{Twisting morphism for $\Po$-algebras}\label{algtwbar}

In this section, we extend the Koszul duality theory for associative algebras to algebras over an operad. We recall the notions already in \cite{GetzlerJones} of \emph{algebraic twisting morphism} and the \emph{bar and the cobar constructions}. However, to describe the Koszul duality theory for $\Po$-algebras, we need to generalize the Koszul complex. A cohomology theory associated to $\Po$-algebras is represented by the \emph{cotangent complex}, that we make explicit thanks to the Koszul duality theory for operads. We show that this cotangent complex generalizes the Koszul complex and we state and prove in this setting the algebraic twisting morphisms fundamental theorem.

We fix an augmented dg operad $\Po$, a coaugmented dg cooperad $\C$ and an operadic twisting morphism $\alpha : \C \rightarrow \Po$. From now on and until the end of the paper, we assume that $\Po(0) = 0$ and $\C(0)=0$.

\subsection{$\Po$-algebra}

A \emph{$\Po$-algebra} is a dg module $A$ endowed with a morphism of dg operads
$$\Po \rightarrow End_{A} := \{ \mathrm{Hom}(A^{\otimes n},\, A)\}_{n\geq 0}.$$
Equivalently, a structure of $\Po$-algebra $(A,\, \gamma_{A})$ is given by a map $\gamma_A : \Po (A) \rightarrow A$ which is compatible with the composition product of the operad $\Po$ and the unit of the operad $\Po$, where
$$\Po (A) := \Po \circ (A,\, 0,\, 0,\, \cdots) = \bigoplus_{n\geq 0} \Po(n) \otimes_{\So_{n}} A^{\otimes n}.$$
Dually, a structure of \emph{$\C$-coalgebra  $(C,\, \Delta_{C})$} is a dg module $C$ endowed with a map
$$\Delta_{C} : C \rightarrow \C (C) := \prod_{n\geq 0} (\C(n) \otimes C^{\otimes n})^{\So_{n}}$$
which is compatible with the coproduct and the counit and where $(-)^{\So_{n}}$ stands for the coinvariants with respect to the diagonal action. We say that the $\C$-coalgebra $C$ is \emph{conilpotent} when the map $\Delta_{C}$ factors through $\oplus_{n \geq 0} (\C(n) \otimes C^{\otimes n})^{\So_{n}}$.

The notation $\otimes_{H}$ stands for the \emph{Hadamard product}: for any $\So$-modules $M$ and $N$, $(M\otimes_{H} N)(n) := M(n) \otimes N(n)$ with the diagonal action of $\So_{n}$. Let ${\Sc^{-1}}$ be the cooperad $\mathrm{End}^{c}_{s\K} := \{\textrm{Hom}((s\K)^{\otimes n},\, s\K)\}_{n \geq 0}$ endowed with a natural action of $\So_{n}$ given by the signature, where $s$ stands for the homological suspension of vector spaces. To the cooperad $\C$, we associate its \emph{operadic homological desuspension} given by the cooperad $\Cs := \Scc \otimes_{H} \C$. A structure of $\Cs$-coalgebra $\Delta_{sC}$ on $sC$ is equivalent to a structure of $\C$-coalgebra $\Delta_{C}$ on $C$ because $\Scc \C (sC) \cong s\C (C)$.

\begin{ex1}
Let $\Po = \As$ be the non-symmetric operad encoding associative algebras. An $\As$-algebra $A$ is an associative algebra without unit. Moreover, a structure of $\Scc \As^{\ash}$-coalgebra is exactly that of a coassociative coalgebra as in Section \ref{casas} but without counit. The category of $\As$-algebras $A$ is equivalent to the category of augmented associative algebra by adding a unit $A_{+} := \K \oplus A$. The Koszul duality theory of $\As$-algebras will be the classical Koszul duality theory of augmented associative algebras.
\end{ex1}

A \emph{weight graded dg module}, or \emph{wdg module} for short, is a chain complex endowed with a weight grading. We say that a \emph{wdg $\Po$-algebra} or \emph{wdg $\C$-coalgebra} $V$ is \emph{connected} when it satisfies $V = V^{(1)} \oplus V^{(2)} \oplus \cdots$. Moreover, we require that the structure maps, as the composite product $\gamma_{A}$, preserve the weight grading.

\subsection{Algebraic twisting morphism}\label{algtwmorph}

From \cite{GetzlerJones}, we recall the definition of twisting morphisms between a coalgebra over a cooperad and an algebra over an operad. We describe the bar and the cobar constructions in this setting.\\

Let $(A,\, \gamma_{A})$ be a $\Po$-algebra and $(C,\, \Delta_{C})$ be a $\C$-coalgebra. Associated to $\varphi \in \mathrm{Hom_{dg\, mod}}(sC,\, A)$, we define the applications
$$\left\{ \begin{array}{l}
\star_{\alpha} (\varphi) : sC \xrightarrow{\Delta_{C}} s\C (C) \xrightarrow{(s^{-1}\alpha) \circ (s\varphi)} \Po (A) \xrightarrow{\gamma_{A}} A\\
\partial(\varphi) := d_{A} \circ \varphi - (-1)^{|\varphi|} \varphi \circ d_{sC}.
\end{array} \right.$$

An \emph{algebraic twisting morphism with respect to $\alpha$} is a map $\varphi : sC \rightarrow A$ of degree $-1$ solution to the Maurer-Cartan equation
$$\partial (\varphi) + \star_{\alpha} (\varphi) = 0.$$
We denote by Tw$_{\alpha}(C, A)$ the set of algebraic twisting morphisms with respect to $\alpha$. In the weight graded setting, we require that the algebraic twisting morphisms preserve the weight grading.

\begin{rem}
We proved in \cite{Milles} that when $\Po$ is binary quadratic operad and $\C = \Poa$ is the Koszul dual cooperad, the aforementioned Maurer-Cartan equation is equal to a Maurer-Cartan equation in a dg Lie algebra.
\end{rem}

To the operadic twisting morphism $\alpha : \C \rightarrow \Po$, one associates a functor
$$\B_{\alpha} : \textsf{dg}\ \Po \textsf{-algebras} \rightarrow \textsf{dg}\ \Cs \textsf{-coalgebras}$$
defined by $\B_{\alpha} A = s\C \circ_{\alpha} A := (s\C (A),\, d_{\alpha} := id_{s} \otimes (d_{\C(A)} + d_{\alpha}^{r}))$, where
$$d_{\alpha}^{r} : \C (A) \xrightarrow{\Delta_{(1)} \circ id_{A}} (\C \circ_{(1)} \C) (A) \xrightarrow{(id_{\C} \circ_{(1)} \alpha) \circ id_{A}} \C \circ \Po (A) \xrightarrow{id_{\C}\circ \gamma_{A}} \C (A).$$
The coproduct defining the $\Scc \C$-coalgebra structure on $s\C(A)$ is given by
$$s\C(A) \xrightarrow{s\Delta_{\C} \circ id_{A}} s\C \circ \C (A) \cong \Cs(s\C (A)).$$

\begin{lem}[\cite{GetzlerJones}]
The map $d_{\alpha}$ is a differential, that is $d_{\alpha}^{2} = 0$ and $\B_{\alpha} A$ is a dg $\Cs$-coalgebra.
\end{lem}

In a similar way, one associates to the operadic twisting morphism $\alpha : \C \rightarrow \Po$ a functor
$$\Omega_{\alpha} : \textsf{dg}\ \Cs \textsf{-coalgebras} \rightarrow \textsf{dg}\ \Po \textsf{-algebras}$$
given on a $\Cs$-coalgebra $sC$ by $\Omega_{\alpha} sC := (\Po(C),\, d_{\alpha} := d_{\Po(C)} - d_{\alpha}^{l})$, where
$$d_{\alpha}^{l} : \Po(C) \xrightarrow{id_{\Po}\circ' \Delta_{C}} \Po \circ (C,\, \C (C)) \xrightarrow{id_{\Po} \circ (id_{C},\, \alpha \circ id_{C})} (\Po \circ_{(1)} \Po) (C) \xrightarrow{\gamma_{(1)} \circ id_{\C}} \Po (C).$$

\begin{lem}[\cite{GetzlerJones}]
The map $d_{\alpha}$ is a differential, that is $d_{\alpha}^{2} = 0$ and $\Omega_{\alpha} sC$ is a $\Po$-algebra.
\end{lem}

Notice that the notation $d_{\alpha}$ stands for different differentials. The differential is given without ambiguity by the context.

\begin{ex1}
Assume that $\Po = \As$ is the operad encoding associative algebras, $\C = \As^{\ash}$ and $\kappa : \As^{\ash} \rightarrow \As$. Let $A$ be an $\As$-algebra, that is an associative algebra.

The bar construction $\B_{\kappa}A$ of the $\As$-algebra $A$ is equal to the classical (reduced) bar construction $\overline{\B} A_{+} := (\overline{T}^{c}(sA),\, d)$ \cite{EilenbergMacLane} of the associative algebra $A_{+} = \K \oplus A$. Moreover, let $C$ be an $\As^{\ash}$-coalgebra, that is $sC$ is a coassociative coalgebra without counit. The cobar construction $\Omega_{\kappa}sC$ is equal to the classical (reduced) cobar construction $\overline{\Omega} sC_{+} := (\overline{T}(C),\, d)$ \cite{Adams} of the coaugmented coassociative coalgebras $sC_{+} = \K \oplus sC$.
\end{ex1}

The bar construction and the cobar construction form the \emph{bar-cobar adjunction}.

\begin{prop}[Proposition $2.18$ of \cite{GetzlerJones}]\label{bijbar}
For every conilpotent $\C$-coalgebra $C$ and every $\Po$-algebra $A$, there is a natural bijection
$$\begin{array}{c}
\mathrm{Hom}_{\mathsf{dg\, \Po} \textsf{-}\mathsf{alg.}}(\Omega_{\alpha}sC,\, A) \cong \mathrm{Tw}_{\alpha}(C,\, A) \cong \mathrm{Hom}_{\mathsf{dg\, \Cs} \textsf{-}\mathsf{cog.}}(sC,\, \B_{\alpha}A)\\
\xymatrix@C=40pt@M=8pt{f_{\varphi} & \ar@{<->}[l] \varphi \ar@{<->}[r] &  g_{\varphi}.}
\end{array}$$
\end{prop}

This adjunction produces two particular morphisms. Consider a $\Po$-algebra $A$ and its bar construction $sC = \B_{\alpha}A$. The morphism of dg $\Cs$-coalgebras $id_{\B_{\alpha}A}$ gives a universal algebraic twisting morphism
$$\pi_{\alpha} : \B_{\alpha}A \cong s\C(A) \epi sA \xrightarrow{s^{-1}} A$$
and the \emph{counit} of the adjunction
$$\varepsilon_{\alpha} : \Omega_{\alpha}\B_{\alpha}A = \Po \circ_{\alpha} \C \circ_{\alpha} A \xrightarrow{id_{\Po} \circ (s\pi_{\alpha})} \Po (A) \xrightarrow{\gamma_{A}} A.$$
Similarly, when $C$ is a dg $\C$-coalgebra and $A = \Omega_{\alpha}sC$, the morphism $id_{\Omega_{\alpha}sC}$ of dg $\Po$-algebras gives a universal algebraic twisting morphism
$$\iota_{\alpha} : sC \xrightarrow{s^{-1}} C \mono \Omega_{\alpha}sC \cong \Po (C)$$
and the \emph{unit} of the adjunction
$$u_{\alpha} : sC \xrightarrow{\Delta{sC}} s\C(C) \xrightarrow{id_{s\C} \circ (s\iota_{\alpha})} \B_{\alpha}\Omega_{\alpha}\, sC.$$
The morphisms $\pi_{\alpha}$ and $\iota_{\alpha}$ are universal in the following meaning.

\begin{lem}[\cite{GetzlerJones}]
With the above notations, any algebraic twisting morphism $\varphi : sC \rightarrow A$ with respect to $\alpha$ factors through the universal algebraic twisting morphisms
$$\xymatrix{& \Omega_{\alpha}sC \ar@{-->}[dr]^{f_{\varphi}} &\\
sC \ar[rr]^{\varphi} \ar[ur]^{\iota_{\alpha}} \ar@{-->}[dr]_{g_{\varphi}} && A\\
& \B_{\alpha}A. \ar[ur]_{\pi_{\alpha}} &}$$
\end{lem}

\begin{pf}
The dashed arrows are just the images of $\varphi$ by the two bijections of Proposition \ref{bijbar}.
$\cqfd$
\end{pf}

We now prove that the bar and the cobar construction behave well in the weight graded setting.

\begin{lem}\label{cobarbarpreserve}
Let $\alpha : \C \rightarrow \Po$ be a Koszul morphism between a wdg connected cooperad $\C$ and a wdg connected operad $\Po$. The cobar construction $\Omega_{\alpha}$ sends quasi-isomorphisms $g : sC \qiso sC'$ between wdg connected $\Scc \C$-coalgebras to quasi-isomorphisms $\Omega_{\alpha} sC \qiso \Omega_{\alpha} sC'$ of $\Po$-algebras.

Similarly, the bar construction $\B_{\alpha}$ sends quasi-isomorphisms $f : A \qiso A'$ between wdg connected $\Po$-algebras to quasi-isomorphisms $\B_{\alpha} A \qiso \B_{\alpha} A'$ of $\Cs$-coalgebras.
\end{lem}

\begin{pf}
We show first the result for the cobar construction. Since $\K$ is a field of characterisic $0$, every dg module is projective. Moreover, by Maschke's theorem, every $\K[\So_{n}]$-module is projective. So the quasi-isomorphism $sC \qiso sC'$ implies the quasi-isomorphism $C \qiso C'$, $C^{\otimes n} \qiso {C'}^{\otimes n}$ and $\Po (n) \otimes_{\So_{n}} C^{\otimes n} \qiso \Po (n) \otimes_{\So_{n}} {C'}^{\otimes n}$. Finally $(\Po (C),\, d_{\Po(C)}) \qiso (\Po (C'),\, d_{\Po(C')})$. We filter the chain complex $\Omega_{\alpha}sC = (\Po (C),\, d_{\Po(C)} - d_{\alpha}^{l})$ by the total weight in $C$
$$F_{p}(\Po (C)) := \bigoplus_{\omega_{1}+ \cdots +\omega_{m}\leq p} \bigoplus_{n \in \mathbb{N}} \Po(n) \otimes_{\So_{n}} C^{(\omega_{1})} \otimes \cdots \otimes C^{(\omega_{n})}.$$
The part $d_{\Po(C)}$ of the differential keeps the total weight in $C$ constant. Since $\C$ and $\Po$ are wdg connected, the twisting morphisms are zero on weight zero and the part $d_{\alpha}^{l}$ of the differential decreases the total weight at least by $1$. The differential respects the filtration. This filtration is exhaustive and bounded below. So we apply the classical theorem of convergence of spectral sequences (Theorem $5.5.1$ of \cite{Weibel}) and we get that the induced spectral sequence converges to the homology of $\Omega_{\alpha}sC$. We consider the same filtration for $C'$. The terms $E^1_{p,\, q}(\Po (C))$ are given by the homology of $(\Po (C),\, d_{\Po(C)})$ and are isomorphic to the terms $E^1_{p,\, q}(\Po (C'))$, that is to the homology of $(\Po (C'),\, d_{\Po(C')})$. Since moreover $g$ is a morphism of dg $\C$-coalgebras, the pages $E^r$, $r\geq 1$ are isomorphic and $\Omega_{\alpha}sC \qiso \Omega_{\alpha}sC'$ is a quasi-isomorphism.

To prove the result for the bar construction, we consider the filtration $F_{p}$ on $\B_{\alpha} A$ given by
$$F_{p}(s\C (A)) := \bigoplus_{\omega \leq p} \C^{(\omega)}(A).$$
The rest of the proof is similar.
$\cqfd$
\end{pf}

\subsection{Cotangent complex}\label{cotangentcomplex}

Operadic Koszul morphisms provide functorial resolutions of $\Po$-alge\-bras. We use these resolutions to make the cotangent complex involved in the André-Quillen cohomology theory of $\Po$-algebras explicit.\\

The \emph{cotangent complex} associated to a $\Po$-algebra $A$ is a (class of) chain complexes which represents the \emph{André-Quillen cohomology theory of $A$ with coefficients in a module}. We make it explicit following \cite{GoerssHopkins} and \cite{Milles}, where the reader can find complete exposition about the André-Quillen cohomology theory and the cotangent complex.

To a resolution of the $\Po$-algebra $A$ is associated a representation of the cotangent complex of $A$.

\begin{prop}[Proposition 10.3.6 of \cite{LodayVallette}]\label{Koszulandresolution}
The operadic twisting morphism $\alpha : \C \rightarrow \Po$ is Koszul, that is $\Po \circ_{\alpha} \C \circ_{\alpha} \Po \qiso \Po$ if and only if the counit of the adjunction $\Omega_{\alpha}\B_{\alpha}A \qiso A$ is a quasi-isomorphism for every $\Po$-algebra $A$.
\end{prop}

As a consequence, we recover the following theorems.

\begin{thm}[Theorem 2.19 of \cite{GetzlerJones}]
For any $\Po$-algebra $A$, there is a quasi-isomorphism
$$\xymatrix{\Omega_{\pi} \B_{\pi}A = \Po \circ_{\pi} \B \Po \circ_{\pi} A \ar@{->>}[r]^{\hspace{1.6cm}\sim} & A.}$$
\end{thm}

\begin{thm}[Theorem 2.25 of \cite{GetzlerJones}]
When the operad $\Po$ is Koszul, there is a smaller resolution of any $\Po$-algebra $A$
$$\xymatrix{\Omega_{\kappa}\B_{\kappa}A = \Po \circ_{\kappa} \Poa \circ_{\kappa} A \ar@{->>}[r]^{\hspace{1.5cm}\sim} & A.}$$
\end{thm}

To an operadic twisting morphism $\alpha : \C \rightarrow \Po$ and to an algebraic twisting morphism $\varphi : sC \rightarrow A$ with respect to $\alpha$, we associate the following coequalizer $A\otimes^{\Po} C$:
$$\xymatrix{\Po \circ (\Po (A),\, C) \ar@<0.5ex>[r]^{c_{0}} \ar@<-0.5ex>[r]_{c_{1}} & \Po \circ (A,\, C) \ar@{->>}[r] & A\otimes^{\Po} C,}$$
where $\left\{ \begin{array}{l}
c_0 : \Po \circ (\Po(A),\, C) \rightarrow (\Po \circ \Po)(A,\, C) \xrightarrow{\gamma (id_{A},\, id_{C})} \Po (A,\, C)\\
c_1 : \Po \circ (\Po(A),\, C) \xrightarrow{id_{\Po} (\gamma_{A},\, id_{C})} \Po(A,\, C)
\end{array} \right.$.

The differential $d_{\varphi} := d_{A\otimes^{\Po} C} - d_{\varphi}^{l}$ on $A\otimes^{\Po} C$ depends on the differentials on $A$, $\Po$ and $C$ and on a twisting term $d_{\varphi}^{l}$, which is the map on the quotient $A\otimes^{\Po} C$ induced by $d_{1}^{l} := \sum d_{1}^{l}(n)$ where
$$\hspace{-2cm} d_{1}^{l}(n) : \Po (A,\, C) \xrightarrow{id_{\Po} (id_{A},\, \Delta_{C}(n))} \Po (A,\, (\C(n)\otimes C^{\otimes n})^{\So_{n}}) \xrightarrow{id_{\Po} (id_{A},\, \alpha \otimes (s\varphi)^{\otimes n-1} \otimes id_{C})}$$
$$\hspace{4cm}\Po (A,\, \Po(n) \otimes A^{\otimes n-1} \otimes C) \rightarrow (\Po \circ \Po) (A,\, C) \xrightarrow{\gamma (id_{A},\, id_{C})} \Po (A,\, C),$$
with $\Delta_{C}(n) : C \xrightarrow{\Delta_{C}} \C(C) \epi (\C(n) \otimes C^{\otimes n})^{\So_{n}}$ (see Section 2 in \cite{Milles} for more details).

When $\Omega_{\alpha} sC = \Po \circ_{\alpha} C \qiso A$ is a resolution of the $\Po$-algebra $A$, the chain complex $A\otimes^{\Po} C$ is a representation of the \emph{cotangent complex}. For example, when $\alpha$ is a Koszul morphism, the universal twisting morphism $\pi_{\alpha} : \B_{\alpha}A \rightarrow A$ gives the functorial resolution $\Omega_{\alpha} \B_{\alpha} A \cong \Po \circ_{\alpha} \C \circ_{\alpha} A \qiso A$ and a representation of the cotangent complex is given by $A\otimes^{\Po} \B_{\alpha}A$.

\begin{rem}
We denote it by $A\otimes^{\Po} \B_{\alpha}A$ instead of $A\otimes^{\Po} s^{-1}\B_{\alpha}A$ to simplify the notation.
\end{rem}

\begin{ex1}
In the case $\Po = \As$, $\C = \As^{\ash}$ and $\alpha = \kappa : \As^{\ash} \rightarrow \As$. The cotangent complex $A \otimes^{\As} \B_{\kappa}A$ is equal to the augmented bar construction $A_{+} \otimes \overline{\B} A_{+} \otimes A_{+}$, where $A_{+} := \K \oplus A$ since $\B_{\kappa} A = \overline{B} A_{+}$ (see \cite{Milles} for more details).
\end{ex1}

\subsection{Algebraic Koszul morphisms}\label{algKosmorph}

Let $\alpha : \C \rightarrow \Po$ be an operadic Koszul morphism and let $\varphi : sC \rightarrow A$ be an algebraic twisting morphism. The associated $\Scc \C$-algebras morphism $g_{\varphi} : sC \rightarrow \B_{\alpha} A$ induces a natural morphism of dg $A$-modules $A\otimes^{\Po} C \rightarrow A\otimes^{\Po} \B_{\alpha}A$. We say that $\varphi$ is an \emph{algebraic Koszul morphism} when the morphism $A\otimes^{\Po} C \rightarrow A\otimes^{\Po} \B_{\alpha}A$ is a quasi-isomorphism. We denote by $\mathrm{Kos}_{\alpha}(C,\, A)$ the set of algebraic Koszul morphisms from $sC$ to $A$.

There are some operads $\Po$ such that the André-Quillen cohomology theory of any $\Po$-algebra $A$ is an Ext-functor over the enveloping algebra of $A$. These operads satisfy the following property:
\begin{equation}
\textrm{There is a quasi-isomorphism $A \otimes^{\Po} \B_{\alpha}A \qiso \Omega_{\Po}(A)$ for any $\Po$-algebra $A$,}\tag{$\star$}
\end{equation}
where $\Omega_{\Po}(A)$ is the module of K\"ahler differentials forms. We refer to \cite{Milles} for the complete study. Hence, in this case, an algebraic twisting morphism $\varphi : sC \rightarrow A$ is an algebraic Koszul morphism if and only if the map $A \otimes^{\Po} C \qiso \Omega_{\Po}(A)$ is a quasi-isomorphism.

\begin{thm}[Algebraic twisting morphisms fundamental theorem]\label{atmft}
Let $\alpha : \C \rightarrow \Po$ be a Koszul morphism between a wdg connected cooperad $\C$ and a wdg connected operad $\Po$. Let $C$ be a wdg connected $\C$-coalgebra and $A$ be a wdg connected $\Po$-algebra. Let $\varphi : sC \rightarrow A$ be an algebraic twisting morphism. The following assertions are equivalent:
\begin{enumerate}
\item the twisting morphism $\varphi$ is an algebraic Koszul morphism, that is $A\otimes^{\Po} C \qiso A\otimes^{\Po} \B_{\alpha}A$;
\item the map of $\Cs$-coalgebras $g_{\varphi} : sC \qiso \B_{\alpha}A$ is a quasi-isomorphism;
\item the map of $\Po$-algebras $f_{\varphi} : \Omega_{\alpha}sC \qiso A$ is a quasi-isomorphism.
\end{enumerate}
Moreover, when $\Po$ satisfies Condition $(\star)$, the previous assertions are equivalent to
\begin{enumerate}
\item[(1')] the natural map $A\otimes^{\Po} C \qiso \Omega_{\Po}(A)$ is a quasi-isomorphism;
\end{enumerate}
\end{thm}

\begin{rem}
When the operad $\Po$ and the $\Po$-algebra $A$ are concentrated in homological degree $0$, the module of K\"ahler differential forms $\Omega_{\Po}(A)$ is concentrated in homological degree $0$ and Condition $(1')$ writes: $A \otimes^{\Po} C$ is acyclic.
\end{rem}

\begin{pf}
We apply the comparison Lemma proved in Appendix \ref{comparisonlemma} to $sC$ and $sC' = \B_{\alpha}A$ to get the equivalence between $(1)$ and $(2)$.

To prove the equivalence $(2) \Leftrightarrow (3)$, we apply Lemma \ref{cobarbarpreserve} to the quasi-isomorphism $sC \qiso \B_{\alpha}A$ to get the quasi-isomorphism $\Omega_{\alpha}sC \qiso \Omega_{\alpha}\B_{\alpha}A$. Since $\alpha$ is a Koszul morphism, $\Omega_{\alpha}\B_{\alpha}A \qiso A$ by Proposition \ref{Koszulandresolution} and we get the implication $(2) \Rightarrow (3)$.

To prove the reverse implication, we apply Lemma \ref{cobarbarpreserve} to the quasi-isomorphism $\Omega_{\alpha} sC \qiso A$ to get the quasi-isomorphism $\B_{\alpha} \Omega_{\alpha}\, sC \qiso \B_{\alpha}A$. Then we just need to prove that $\B_{\alpha} \Omega_{\alpha} sC = s\C \circ_{\alpha} \Po \circ_{\alpha} C \qiso sC$ for each $\C$-coalgebra $C$, provided that $\C \circ_{\alpha} \Po \qiso I$ (by Theorem \ref{otmft}). To prove this, we endow $\C \circ_{\alpha} \Po \circ_{\alpha} C$ with a filtration $F_{p}$ given by
$$F_{p}(\C \circ_{\alpha} \Po \circ_{\alpha} C) := \bigoplus_{\omega \leq p} \C \circ_{\alpha} \Po \circ_{\alpha} \underbrace{C}_{(\omega)}.$$
The differential $d_{\C \circ \Po \circ C} + d_{\alpha}^{r} \circ id_{C} - id_{\C} \circ d_{\alpha}^{l}$ satisfies
$$\left\{ \begin{array}{lcl}
d_{\C \circ \Po \circ C} & : & F_{p} \rightarrow F_{p}\\
d_{\alpha}^{r} \circ id_{C} & : & F_{p} \rightarrow F_{p}\\
id_{\C} \circ d_{\alpha}^{l} & : & F_{p} \rightarrow F_{p-1}.
\end{array} \right.$$
So the filtration is a filtration of chain complexes. Moreover, it is bounded below and exhaustive so the classical theorem of convergence of spectral sequences (Theorem 5.5.1 of \cite{Weibel}) gives that the induced spectral sequence $E_{p,\, q}^{\bullet}$ converges to the homology of $\C \circ_{\alpha} \Po \circ_{\alpha} C$. We endow similarly $C$ with the filtration by the weight grading and we get an isomorphism between the $E^{1}$-pages, provided that $E_{p,\, q}^{1} \cong \mathrm{H}_{\bullet}(\C \circ_{\alpha} \Po) \circ \mathrm{H}_{\bullet}(C) \cong \mathrm{H}_{\bullet}(C)$ by the K\"unneth formula. The convergence of the spectral sequences concludes the proof.
$\cqfd$
\end{pf}

\begin{ex1}
When $\Po = \As$, we recover Theorem \ref{tmft}, provided that the category of $\As$-algebras is equivalent to the category of augmented associative algebras and that a conipoltent $\Scc \As$-coalgebra is exactly a conilpotent coaugmented coassociative coalgebra.
\end{ex1}

\section{Koszul duality theory for algebra over an operad}

We proved in the previous section the algebraic twisting morphisms fundamental theorem. In this section, we define the notion of \emph{monogene $\Po$-algebra} and we associate to such a $\Po$-algebra $A$ its Koszul dual $\Scc \Poa$-coalgebra , which is a good candidate for the algebraic twisting morphisms fundamental theorem. We show the link with the Koszul dual $\Po^{!}$-algebra defined in \cite{GinzburgKapranov}. We give the Koszul criterion for $\Po$-algebra and the definition of a Koszul $\Po$-algebra. Thus we obtain a criterion to prove that we have a ``small'' resolution of $A$. This generalizes the Koszul duality theory for quadratic associative algebras \cite{Priddy} and conceptually explains the form of the Koszul complex by the fact that it is a representation of the cotangent complex when the $\Po$-algebra is Koszul.

\subsection{Monogene $\Po$-algebra}

The notion of quadratic algebra over a quadratic operad appears in \cite{GinzburgKapranov}. We extend this definition to \emph{monogene} $\Po$-algebra and we give the definition of \emph{monogene} $\C$-coalgebra. We use the word monogene to express the fact that the relations in the $\Po$-algebra are linearly generated by the generators of the operad.\\

Let $(E,\, R)$ be an operadic quadratic data (see Section \ref{Koszuloperad}) and $\Po := \F(E)/(R)$ its associated quadratic operad. A \emph{monogene data $(V,\, S)$ of the operadic quadratic data $(E,\, R)$} is a graded vector space $V$ and a subspace $S \subseteq E(V)$. We associate to this monogene data the \emph{monogene $\Po$-algebra}
$$A(V,\, S) := \Po(V)/(S).$$

The operad $\Po$ is weight graded and we endow $\Po(V)$ with a weight grading equal to the weight grading in $\Po$. The space $S$ is homogeneous for this grading since the weight of $E(V)$ is equal to $1$, so the $\Po$-algebra $A(V,\, S)$ is weight graded.

\begin{prop}
When $\Po = \Po(E,\, R)$ is a binary quadratic operad, that is $E = E(2) = \Po(2)$, we have $\Po(0) = 0$ and $\Po(1) = \K$ and explicitly
$$A(V,\, S) = \bigoplus_{n\geq 0}A(V,\, S)^{(n)} = V \oplus (E \otimes_{\So_{2}} V^{\otimes 2})/S \oplus \cdots \oplus (\Po (n) \otimes_{\So_{n}} V^{\otimes n})/ \mathrm{Im}(\psi_{n}) \oplus \cdots,$$
where $\psi_{n}$ is the composite
$$\bigoplus_{j=0}^{n-1} \Po (n-1) \otimes_{\So_{n-1}}(V^{\otimes j} \otimes S \otimes V^{\otimes n-2-j}) \rightarrow (\Po \circ_{(1)} E) (n) \otimes_{\So_{n}} V^{\otimes n} \mono \Po \circ \Po (V)\xrightarrow{\gamma \circ id_{V}} \Po (V).$$
In this case, the monogene $\Po$-algebra is a \emph{quadratic $\Po$-algebra} as defined in \cite{GinzburgKapranov}.
\end{prop}

\begin{pf}
The image of the map $\psi_{n}$ is exactly the definition of the ideal of $\Po(V)$ generated by $S \subseteq E(V)$.
$\cqfd$
\end{pf}

\begin{ex1}
When $\Po = \As$, we recover the notion of quadratic associative algebra $T(V)/(S)$ and the weight grading is given by the number of elements in $V$ minus $1$.
\end{ex1}

Dually, let $\C = \C(E,\, R)$ be the cooperad associated to the operadic quadratic data $(E,\, R)$. The \emph{monogene $\C$-coalgebra} associated to the monogene data $(V,\, S)$ of the operadic quadratic data $(E,\, R)$ is the $\C$-coalgebra $C(V,\, S)$ which is universal among the sub-$\C$-coalgebras $C$ such that the composite
$$C \mono \C(V) \epi E(V)/S$$
is equal to $0$. The word ``universal'' means that for any such $\C$-coalgebra $C$, there exists a unique morphism of $\C$-coalgebras $C \rightarrow C(V,\, S)$ such that the following diagram commutes
$$\xymatrix@R=12pt@M=8pt@W=12pt@H=1pt{C(V,\, S) \ar@{>->}[r] & \C (V).\\
C \ar@{>->}[ur] \ar[u] & }$$
The cooperad $\C$ is weight graded, so the $\C$-coalgebra $C(V,\, S)$ is weight graded by the weight on $\C$.

\begin{prop}
When $\C = \C(E,\, R)$ is a binary quadratic cooperad, that is $E = E(2) = \C(2)$, we have $\C(0) = 0$ and $\C(1) = \K$. Dually to the algebra case, we have explicitly
$$C(V,\, S) = \bigoplus_{n\geq 0}C(V,\, S)^{(n)} = V \oplus S \oplus \cdots \oplus \mathrm{Ker}(\phi_{n}) \oplus \cdots,$$
where $\phi_{n}$ is the composite\\
$\C(V) \xrightarrow{\Delta \circ id_{V}} \C \circ \C (V) \epi ((\C \circ_{(1)} E)(n) \otimes V^{\otimes n})^{\So_{n}} \rightarrow$
$$\hspace{2.7cm} ((\C \circ_{(1)} E)(n) \otimes V^{\otimes n})^{\So_{n}}/\bigcap_{j=0}^{n-1} (\C(n-1) \otimes V^{\otimes j} \otimes S \otimes V^{n-2-j})^{\So_{n}}.$$
In this case, the monogene $\C$-coalgebra is a \emph{quadratic $\C$-coalgebra}.
\end{prop}

\begin{pf}
We dualize the map $\psi_{n}$ of the previous proposition to get the map $\phi_{n}$ and the notion of ``coideal'' of $\C(V)$ ``cogenerated'' by $S$.
$\cqfd$ 
\end{pf}

\subsection{Koszul dual coalgebra}

We define the \emph{Koszul dual $\Pos$-coalgebra}, the \emph{Koszul dual $\Po^!$-algebra} of a monogene $\Po$-algebra and the corresponding algebraic twisting morphism. When the operad $\Po$ is binary and finitely generated, we recover the definition of the Koszul dual $\Po^{!}$-algebra of \cite{GinzburgKapranov}. We show that the Koszul dual $\Pos$-coalgebra associated to a monogene $\Po$-algebra is the zero homology group for a certain degree of the bar construction of this $\Po$-algebra.\\

To an operadic quadratic data $(E,\, R)$, we associate the operad $\Po := \Po(E,\, R)$, the Koszul dual cooperad $\Poa := \C(sE,\, s^{2}R)$ and its homological desuspension $\Pos := \Scc \otimes_{H} \Poa$. We assume that $\Po$ is a Koszul operad (cf. \ref{Koszuloperad}).

Let $V$ be a vector space and let $S$ be a subspace of $E(V)$. We have $sS \subset sE(V) \mono \Poa(V)$. The \emph{Koszul dual $\Pos$-coalgebra} of the monogene $\Po$-algebra $A(V,\, S)$ is the monogene $\Pos$-coalgebra
$$A^{\ash} := sC(V,\, sS).$$

We define the morphism $\varkappa : A^{\ash} \rightarrow A$ by the linear map of degree $-1$:
$$A^{\ash} := sC(V,\, sS) \epi sV \xrightarrow{s^{-1}} V \mono A(V,\, S) = A.$$

\begin{lem}\label{kappaKos}
We have $\star_{\kappa}(\varkappa) = 0$ and therefore $\varkappa$ is an algebraic twisting morphism.
\end{lem}

\begin{pf}
Due to the definition of $\varkappa$, the term $\star_{\kappa}(\varkappa)$ is $0$ everywhere except maybe on $s^{2}S \subset s^{2}E(V)$ where it is equal to
$$s^2 S \xrightarrow{\Delta_{A^{\ash}}} s^{2}E(V) \xrightarrow{(s^{-1}\kappa) \circ (s\varkappa)} E(V) \xrightarrow{\gamma_{A}} E(V)/S.$$
This last map is $0$ by definition.
$\cqfd$
\end{pf}

\begin{ex1}
When $\Po = \As$, up to adding a unit, we recover the Koszul dual coalgebra defined in \cite{Priddy}.
\end{ex1}

Let $\C = \oplus_{n\geq 0} \C^{(n)}$ be a weight graded cooperad. We define the \emph{weight-graded linear dual} of $\C$ by $\C^* := \oplus_{n\geq 0} {\C^{(n)}}^* = \oplus_{n \geq 0} \mathrm{Hom}_{\K}(\C^{(n)},\, \K)$. The weight-graded linear dual of $\Pos$ is denoted by
$$\Po^! := (\Pos)^* = (\Scc \otimes_{H} \Poa)^*$$
and is a weight graded operad. When $E$ is a finite dimensional $\So_{2}$-module, it corresponds to the operad defined in \cite{GinzburgKapranov} by $\Po^! = \F (E^{\vee})/(R^{\bot})$ (see Theorem 7.6.5 in \cite{LodayVallette} for a proof of this fact). In the coalgebra case, let $C = \oplus_{n\geq 0} C^{(n)}$ be a weight graded $\Pos$-coalgebra. We define the \emph{weight-graded linear dual} of $C$ by $C^* := \oplus_{n\geq 0} {C^{(n)}}^*$. The \emph{Koszul dual $\Po^!$-algebra $A^!$} is the following $\Po^!$-algebra
$$A^! := \left(\bigoplus_{n \geq 0} {s^{-n}A^{\ash}}^{(n)}\right)^* = \bigoplus_{n \geq 0} s^{n}\left({A^{\ash}}^{(n)}\right)^{*}.$$

\begin{prop}
Let $\Po$ be a finitely generated binary Koszul operad. Let $(V,\, S)$ be a finitely generated monogene data. The Koszul dual $\Po^!$-algebra $A^!$ of $A = A(V,\,S)$ is equal to the monogene $\Po^!$-algebra $A^! = A(V^*,\, R^{\bot})$ defined in \cite{GinzburgKapranov} by $V^{*} := \mathrm{Hom}(V,\, \K)$ and $S^{\bot}$ is the annihilator of $S$ for the natural pairing $\langle -,- \rangle : E^{*}(V^{*}) \otimes E(V) \rightarrow \K$.
\end{prop}

\begin{pf}
Since $E$ and $V$ are finite dimensional, we remark that
$$\left(s \bigoplus_{n \geq 0} s^{-n}{\Poa}^{(n)}(V)\right)^{*} = \Po^! (V^*).$$
After a weight graded desuspension, we linearly dualize the exact short sequence
$$0 \rightarrow A^{\ash} \rightarrow s\Poa(V) \rightarrow s^{2}E(V)/s^{2}S \rightarrow 0$$
satisfied by $A^{\ash}$. We get the exact sequence
$$0 \leftarrow A^! \leftarrow \Po^! (V^*) \leftarrow S^{\bot} \leftarrow 0,$$
where the orthogonal space $S^{\bot}$ is the annihilator $S$ for the natural pairing $\langle -,- \rangle : E^{*}(V^{*}) \otimes E(V) \rightarrow \K$. Since $A^{\ash}$ is universal for the first exact sequence, the dual $A^!$ is universal for the second one and is equal to $\Po^{!}(V^{*})/(S^{\bot})$.
$\cqfd$
\end{pf}

Recall that the operad $\Po$ is weight graded by the number of elements in $E$. It induces a weight grading on $A$ (the weight of $V$ being equal to $0$). Hence, we endow the bar construction $\B_{\kappa}A$ with a non-negative \emph{weight-degree} induced only by the weight grading on $A$. We denote it by $\B_{\kappa}^{\omega}A$.

The internal differential on $A$ is $0$ since the differential on $\Po$ and on $V$ are $0$. Thus the differential on $\B_{\kappa}A$ reduces to $d_{\kappa} = id_{s} \otimes d_{\kappa}^{r}$ (defined in Section \ref{algtwmorph}). The differential $d_{\kappa}$ raises the weight-degree by $1$ and we get a cochain complex with respect to this degree. The elements of weight-degree $0$ in $\B_{\kappa}A = s\Poa (A)$ are given by $s\Poa(V)$, then $A^{\ash} \mono \B_{\kappa}^{0}A = s\Poa(V)$.

\begin{prop}
Let $(E,\, R)$ be an operadic quadratic data. Let $(V,\, S)$ be a monogene data. The natural $\Pos$-coalgebras inclusion $g_{\varkappa} : A^{\ash} = sC(V,\, sS) \mono \B_{\kappa}A = \B_{\kappa} A(V,\, S)$ induces an isomorphism of graded $\Pos$-coalgebras
$$\xymatrix@M=8pt{g_{\varkappa} : A^{\ash} \ar@{>->}[r]^{\hspace{-.4cm} \cong} & \mathrm{H}^{0}(\B_{\kappa}^{\bullet}A).}$$
\end{prop}

\begin{pf}
Since there is no element in negative weight-degree, we just need to prove that the inclusion $g_{\varkappa}$ is exactly the kernel of the differential ${id_{s} \otimes d_{\kappa}^{r}}_{|s\Poa(V)}$. The image of $g_{\varkappa}$ lives in weight-degree $0$. Moreover, the morphism $g_{\varkappa}$ commutes with the differentials $d_{A^{\ash}} = 0$ and $d_{\kappa}$ so $d_{\kappa} \circ g_{\varkappa} = g_{\varkappa} \circ d_{A^{\ash}} = 0$ and $A^{\ash} \mono \mathrm{Ker}({d_{\kappa}}_{|weight = 0}) = \mathrm{Ker}({id_{s} \otimes d_{\kappa}^{r}}_{|s\Poa(V)}) =: K$. Since $K$ is the kernel of a $\Pos$-coalgebras morphism, it is a $\Pos$-coalgebra. It is easy to see that the composition $K \mono s\Poa(V) \epi s^{2}E(V)/s^{2}S$ is equal to $0$ since the differential in weight $1$ is the quotient map $s^{2}E(V) \rightarrow s(E(V)/S)$. Due to the universal property of $A^{\ash}$ and since $A^{\ash} \mono K$ is a monomorphism, we get that $A^{\ash} = K = \mathrm{Ker}({id_{s} \otimes d_{\kappa}^{r}}_{|s\Poa(V)})$.
$\cqfd$
\end{pf}

\subsection{Koszul criterion and Koszul $\Po$-algebra}

The previous proposition shows that the Koszul dual $\Pos$-coalgebra is a good candidate to replace the bar construction in the cotangent complex. We state the Koszul criterion shows that it is the case when the algebraic twisting morphism $\varkappa$ is Koszul. We define the notion of \emph{Koszul $\Po$-algebra}.\\

Let $(E,\, R)$ be an operadic quadratic data, $\Po := \Po(E,\, R)$, $\Poa := \C(sE,\, s^{2}R)$ and $\kappa : \Poa \rightarrow \Po$. Let $(V,\, S)$ be a monogene data, $A := A(V,\, S) = \Po(V)/(S)$ the associated monogene $\Po$-algebra and $A^{\ash} := sC(V,\, sS)$ the Koszul dual $\Pos$-coalgebra.

When $\Po$ is a Koszul operad, the bar-cobar construction $\Omega_{\kappa} \B_{\kappa}A$ of $A$ is a cofibrant resolution of $A$. To simplify this resolution, we can replace $\B_{\kappa}A$ by $A^{\ash} \cong \mathrm{H}^{0}(\B_{\kappa}^{\bullet} A)$. This works when $\mathrm{H}^{\bullet}(\B_{\kappa}^{\bullet} A) = \mathrm{H}^{0}(\B_{\kappa}^{\bullet} A)$. The following Koszul criterion shows that it is the case if and only if the algebraic twisting morphism $\varkappa : A^{\ash} \rightarrow A$ is Koszul.

We apply the algebraic twisting morphism fundamental theorem of the previous section (Theorem \ref{atmft}) to get the following theorems, which are the main theorems of Koszul duality for $\Po$-algebras.

\begin{thm}[Koszul criterion]\label{algKoszulcriterion}
Let $(E,\, R)$ be an operadic quadratic data such that $\Po = \Po (E,\, R)$ is a Koszul operad. Let $(V,\, S)$ be a monogene data associated to $(E,\, R)$. The following assertions are equivalent:
\begin{enumerate}
\item the twisting morphism $\kappa$ is an algebraic Koszul morphism, that is
$$A \otimes^{\Po} A^{\ash} \qiso A \otimes^{\Po} \B_{\kappa}A$$
is a quasi-isomorphism;
\item the inclusion $g_{\varkappa} : A^{\ash} \mono \B_{\kappa}A$ is a quasi-isomorphism;
\item the projection $f_{\varkappa} : \Omega_{\kappa}A^{\ash} \epi A$ is a quasi-isomorphism.
\end{enumerate}
Moreover, when $\Po$ satisfies Condition $(\star)$, the previous assertions are equivalent to
\begin{enumerate}
\item[(1')] the natural map $A\otimes^{\Po} A^{\ash} \qiso \Omega_{\Po}(A)$ is a quasi-isomorphism.
\end{enumerate}
When these assertions hold, the cobar construction on $A^{\ash}$ gives a cofibrant resolution of the $\Po$-algebra $A$ (a minimal resolution when $\Po(0) = 0$ and $\Po(1) = \K$).
A monogene $\Po$-algebra $A$ is called \emph{Koszul} when it satisfies the equivalent properties of this theorem.
\end{thm}

\begin{pf}
We apply Lemma \ref{kappaKos} and Theorem \ref{atmft} to $\C := \Poa$, $A := A(V,\, S)$, $C := A^{\ash}$ and $\varphi := \varkappa$ since the weight assumptions are satisfied.

When the assertions of the theorem hold, $\Omega_{\kappa}A^{\ash}$ is a cofibrant resolution of $A$ since it is a resolution of $A$ and since $\Omega_{\kappa}A^{\ash}$ is a quasi-free $\Po$-algebra on the connected weight graded $\As$-coalgebra $s^{-1}A^{\ash}$. Moreover, when the operad satisfies $\Po(0)=0$ and $\Po(1)=0$, the differential satisfies
$$d_{\Omega_{\kappa}A}(s^{-1}A^{\ash}) = -d_{\kappa}^{l}(s^{-1}A^{\ash}) \subset \oplus_{n\geq 2}\Po(n) \otimes_{\So_{n}} (s^{-1}A^{\ash})^{\otimes n}$$
by construction and $\Omega_{\kappa}A$ is a minimal resolution of $A$.
$\cqfd$
\end{pf}

\begin{rem}
We recall that when the operad $\Po$ and the $\Po$-algebra $A$ are concentrated in homological degree $0$, the module of K\"ahler differential forms $\Omega_{\Po}(A)$ is concentrated in homological degree $0$ and the condition $(1')$ writes simply: $A \otimes^{\Po} A^{\ash}$ is acyclic.
\end{rem}

The chain complex $A\otimes^{\Po} A^{\ash}$ is called the \emph{Koszul complex}. Thus, when the algebraic twisting morphism $\varkappa : A^{\ash} \rightarrow A$ is an algebraic Koszul morphism, the Koszul complex $A\otimes^{\Po} A^{\ash}$ is a representation of the cotangent complex, representing the cohomology theory of the $\Po$-algebra $A$.

\begin{ex1}
Assume that $\Po = \As$, $\C = \As^{\ash}$ and $\kappa : \As^{\ash} \rightarrow \As$ is the operadic twisting morphism between them. The operad $\As$ satisfies the Condition $(\star)$ of Section \ref{algKosmorph}, so an $\As$-algebra $A$ is Koszul if and only if for $A_{+} := \K \oplus A$, the $A$-bimodules morphism $A_{+} \otimes s^{-1}\overline{(A_{+})^{\ash}} \otimes A_{+} \cong A \otimes^{\As} A^{\ash} \qiso \Omega_{\As}A \cong A \otimes A_{+}$ is a quasi-isomorphism. This is equivalent to $A_{+} \otimes (A_{+})^{\ash} \otimes A_{+} \cong s(A \otimes^{\As} A^{\ash}) \oplus (A_{+} \otimes A_{+}) \qiso A_{+}$ is a quasi-isomorphism. This last quasi-isomorphism is the classical definition for a quadratic augmented associative algebra $A_{+}$ to be Koszul (see Theorem \ref{Koszulcriterion} and the definition after).\\
\end{ex1}

For an associative algebra $A$, we already know that $A$ is Koszul if and only if $A^{!}$ is Koszul. As an application of the Koszul criterion theorem, we obtain the same result for $\Po$-algebras.

\begin{thm}\label{AAdual}
Let $(E,\, R)$ be a finitely generated binary operadic quadratic data such that $\Po = \Po (E,\, R)$ is a Koszul operad. Let $(V,\, S)$ be a finitely generated monogene data associated to $(E,\, R)$. The $\Po$-algebra $A := A(V,\, S)$ is Koszul if and only if the $\Po^{!}$-algebra $A^{!} = \Po^{!}(V^{*},\, S^{\bot})$ is Koszul.
\end{thm}

\begin{pf}
The Koszul criterion Theorem \ref{algKoszulcriterion} (2) implies that $A$ is a Koszul $\Po$-algebra if and only if $A^{\ash} \qiso \B_{\kappa} A = \Po^{\ash}(A)$ is a quasi-isomorphism. We consider the $n$-desuspension of the weight-degree $n$ part: it is a chain complex morphism (and a quasi-isomorphism of chain complexes)
$$\bigoplus_{n \geq 0} s^{-n}{A^{\ash}}^{(n)} \rightarrow \bigoplus_{n \geq 0} s^{-n}(\Poa (A))^{(n)}$$
since the differentials preserve the weight grading. We linearly dualize the quasi-isomorphisms $s^{-n}{A^{\ash}}^{(n)} \rightarrow s^{-n}(\Poa (A))^{(n)}$ to get the quasi-isomorphisms $\left(s^{-n}{A^{\ash}}^{(n)}\right)^{*} \leftarrow \left(s^{-n}(\Poa (A))^{(n)}\right)^{*}$ since $\K$ is a field. The sum of these quasi-isomorphisms gives the quasi-isomorphism $A^{!} \xleftarrow{\sim} \Po^{!}(A^{*}) = \Po^{!}({{A^{!}}^{!}}^{*}) = \Po^{!}({A^{!}}^{\ash})$. By the Koszul criterion Theorem \ref{algKoszulcriterion} (3), this implies that $A^{!}$ is a Koszul $\Po^{!}$-algebra.
$\cqfd$
\end{pf}

\section{Links and applications}

In this section, we give examples of applications of the present Koszul duality theory. We have seen the case of associative algebras, we describe now the case of commutative algebras and the case of Lie algebras. We also recover the case of modules.

\subsection{The case of commutative algebras}

\subsubsection{\bf Commutative algebras and Lie coalgebras}

Let $\Po = \Com$ be the operad encoding (non necessarily unital) associative and commutative algebras. A $\Com$-algebra structure on $A$ is equivalent to commutative and associative algebra structure on $A$ given by a commutative product $\gamma_{A} : A^{\otimes 2} \rightarrow A$ satisfying the associativity relation.

\smallskip

The Koszul dual cooperad of $\Com$ is $\Com^{\ash}$, that is the suspension of the cooperad $\Lie^{c}$, which encodes Lie coalgebras. A $\Lie^{c}$-coalgebra structure on $sC$ is equivalent to a Lie coalgebra structure on $sC$ given by an anti-commutative coproduct $\Delta_{sC} : sC \rightarrow s(\Com^{\ash}(2) \otimes C^{\otimes 2})^{\So_{2}} \cong \{c_{1} \otimes c_{2} + (-1)^{|c_{1}||c_{2}|} c_{2} \otimes c_{1};\, c_{1},\, c_{2} \in sC\}$ which satisfies the coJacobi relation (see \cite{LodayVallette} for more details). We remark that $\mathrm{Im} \Delta_{sC} \subset sC \otimes sC$. When $C$ is finitely generated and concentrated in degree $0$, we have $s(\Com^{\ash}(2) \otimes C^{\otimes 2})^{\So_{2}} \cong \Lambda^{2} C$ where $|\Lambda^{2}| = 2$.

\subsubsection{\bf The cotangent complex}

Let $A$ be a $\Com$-algebra, let $sC$ be a Lie coalgebra and let $\varphi : sC \rightarrow A$ be an algebraic twisting morphism. The twisted tensor product $A\otimes^{\Com} C$ is given by
$$A_{+} \otimes_{\varphi} C := (A_{+} \otimes C, d_{\varphi} := d_{A_{+} \otimes C} - d_{\varphi}^{l}),$$
where $A_{+} := \K \oplus A$ is the augmented algebra of $A$. The differential $d_{A_{+} \otimes C}$ is equal to $d_{A_{+}} \otimes id_{C} + id_{A_{+}} \otimes d_{C}$ and the twisting differential $d_{\varphi}^{l}$ is given by
$$A_{+} \otimes C \xrightarrow{id_{A_{+}} \otimes s^{-1}\Delta_{sC}} A_{+} \otimes sC \otimes C \xrightarrow{id_{A_{+}} \otimes \varphi \otimes id_{C}} A_{+} \otimes A \otimes C \xrightarrow{\tilde{\gamma}_{A} \otimes id_{C}} A_{+} \otimes C,$$
where $\tilde{\gamma}_{A} : A_{+} \otimes A \cong A \oplus A \otimes A \xrightarrow{id_{A} + \gamma_{A}} A \mono A_{+}$. We remark that this construction is close to the twisted tensor product of an associative algebra and a coassociative coalgebra.

\smallskip

The algebraic twisting morphism $\varphi$ is Koszul when the $A$-modules morphism $A_{+} \otimes_{\varphi} C \qiso A_{+} \otimes_{\pi_{\kappa}} s^{-1}\B_{\kappa} A \cong A_{+} \otimes_{\pi_{\kappa}} s^{-1}\Lie^{c}(sA)$, where $\Lie^{c}(sA)$ is the cofree Lie coalgebra on $sA$, is a quasi-isomorphism.

\smallskip

Under certain assumptions as the smoothness or the regularity of $A$ \cite{Quillen}, the cotangent complex is quasi-isomorphic to the $A$-module of K\"ahler differential forms $\Omega^{1}(A)$, therefore $\varphi$ is a Koszul morphism when the $A$-modules morphism $A_{+} \otimes_{\varphi} C \qiso \Omega^{1}(A)$ is a quasi-isomorphism.

\subsubsection{\bf Quadratic setting}

The operad $\Com$ is binary and therefore a monogene $\Com$-algebra is a quadratic commutative algebra: $A := \mathcal{S}(V)/(S)$, where $\mathcal{S}(V) := \Com(V)$ is the symmetric algebra and where $S \subseteq \mathcal{S}(V)^{(2)} = V \otimes_{\So_{2}} V =: V^{\odot 2}$. The Koszul dual Lie coalgebra $A^{\ash}$ is the quadratic Lie coalgebra $sC(V,\, sS)$. When $V$ is finite dimensional, we desuspend and linearly dualize $A^{\ash}$ to get $A^{!} = \Lie(V^{*})/(R^{\bot})$, where $\Lie(V^{*})$ is the free Lie algebra on $V^{*}$. Provided the fact that $\varkappa : A^{\ash} \rightarrow A$ is a Koszul morphism, the Koszul criterion gives resolutions $A^{\ash} \qiso \B_{\kappa} A$ and $\Omega_{\kappa} A^{\ash} \qiso A$.

\subsubsection{\bf Applications to rational homotopy theory}

In this paper, we work in the homological setting. However, we can reverse the arrows and work in the cohomological setting. When the algebras $A$ are generated by a finitely generated vector space $V$, the theory of this paper still works for commutative algebras. In this example, this condition is always verified.\\

The latter quasi-isomorphism $\Omega_{\kappa} A^{\ash} = \Lambda A^{\ash} \qiso A$ is a quadratic model in rational homotopy theory \cite{Sullivan}. Assume that $X$ is a formal simply connected space and that the cohomology ring $A = \mathrm{H}^{\bullet}(X,\, \mathbb{Q})$ is a finitely generated quadratic algebra $\mathcal{S}(V)/(S)$, where $V$ is homological graded. When $A$ forms a Koszul algebra, the Koszul dual algebra $A^{!} = \Lie (V^{*})/(S^{\bot})$ is equal to the rational homotopy groups of $X$. Therefore the Koszul duality theory generates from the quadratic data all the syzygies; we do not have to compute the syzygies one by one. The Lie algebra structure on $A^{!}$ is the Whitehead Lie bracket. By Theorem \ref{AAdual}, we can either prove that $A$ is Koszul or that $A^{\ash}$ is Koszul. Moreover, the Koszul criterion \ref{algKoszulcriterion} (1) provides a new way to prove these conditions.

\smallskip

The complement of a complex hyperplane arrangement is always a formal space. However, it is not necessarily simply connected. To a complex hyperplane arrangement $\mathcal{A}$, one associates the \emph{Orlik-Solomon algebra} $A := A(\mathcal{A})$. This algebra is naturally isomorphic to the cohomology groups of the complement $X$ of the hyperplane arrangement $\mathcal{A}$, that is $A = A(\mathcal{A}) \cong \mathrm{H}^{\bullet}(X,\, \mathbb{Q})$. The conditions on $\mathcal{A}$ for $A$ to be quadratic are studied in \cite{Yuzvinsky}. When the Orlik-Solomon algebra $A$ is quadratic, the Koszul dual algebra $A^{!}$ is the \emph{holonomy Lie algebra} defined by Kohno \cite{Kohno, Kohno2}. When $A$ is a Koszul algebra, the holonomy Lie algebra $A^{!}$ computes the $n$-homotopy groups of the $\mathbb{Q}$-completion of the space $X$ for $n \geq 2$.

\subsubsection{\bf Relationship with the Koszul duality theory of associative algebras}

The free commutative algebra $\mathcal{S}(V)$ is equal to the quadratic associative algebra $T(V)/\langle v\otimes w - w \otimes v \rangle$, where $T(V)$ is the tensor algebra or the free associative algebra on $V$. Thus, there is a weight preserving projection $p : T(V) \epi \mathcal{S}(V)$ and a functor
$$\begin{array}{lcl}
\textsf{quad. comm. alg.} & \rightarrow & \textsf{quad. assoc. alg.}\\
A := \mathcal{S}(V)/(S) & \mapsto & A_{as} := T(V)/p^{-1}(S).
\end{array}$$

We emphasize the fact that the Koszul complexes associated to $A$ and $A_{as}$ are distinct. The Koszul dual coalgebra $A^{\ash}$ is a Lie coalgebra whereas the Koszul dual coalgebra ${A_{as}}^{\ash}$ is a coassociative coalgebra. However, the Koszul dual algebras are linked by the equality ${A_{as}}^{!} = U(A^{!})$, where $U(A^{!})$ is the enveloping algebra of the Lie algebra $A^{!}$ (see \cite{GinzburgKapranov} for example):
$$\xymatrix@R=20pt{\As \textsf{-alg.} : & A_{as} \ar@{<->}[r]^{\hspace{-.7cm} !} \ar@{<->}[d] & {A_{as}}^{!} = U(A^{!}) \ar@{<->}[d] & : \As \textsf{-alg}\\
\Com \textsf{-alg} : & A \ar@{<->}[r]^{!} & A^{!} & \hspace{.1cm}: \Lie \textsf{-alg}.}$$

A priori, the enveloping algebra of a Lie algebra has quadratic and linear relations. However, when the Lie algebra is a homogeneous quadratic Lie algebra, the enveloping algebra admits a homogeneous quadratic presentation as an associative algebra. In \cite{PapadimaYuzvinsky}, the authors proved the following theorem:
\begin{thm}[Proposition 4.4 of \cite{PapadimaYuzvinsky}]
The quadratic associative algebra $A_{as}$ is Koszul if and only if the quadratic commutative algebra $A$ is Koszul.
\end{thm}

They also proved that for a formal space $X$, the algebra $A = \mathrm{H}^{\bullet}(X,\, \mathbb{Q})$, or $A_{as}$, is Koszul if and only if the space $X$ is a rational $K(\pi,\, 1)$. We emphasize however the fact that the definition for an algebra $A$ to be Koszul in \cite{PapadimaYuzvinsky} is slightly different that the one in this paper. They require the generators of $A$ to be in degree $1$, assumptions which is not required in the present paper.

\subsection{The case of Lie algebras}

Let $\Po = \Lie$ be the operad encoding Lie algebras. The Koszul dual cooperad of $\Lie$ is $\Lie^{\ash}$, that is the suspension of the cooperad $\Com^{c}$ encoding co-commutative coalgebras.

\subsubsection{\bf The twisted tensor product and the cotangent complex}

Let $\g$ be a Lie algebra, let $sC$ be a commutative coalgebra and let $\varphi : sC \rightarrow \g$ be a twisting morphism. The twisted tensor product $\g \otimes^{\Lie} C$ is given by $U(\g) \otimes_{\varphi} C := (U(\g) \otimes C,\, d_{\varphi})$ where $U(\g)$ is the enveloping algebra of the Lie algebra $\g$ and where the differential $d_{\varphi}$ is obtained in the same way as the case of commutative algebras.

\smallskip

The cotangent complex is the Chevalley-Eilenberg complex $U(\g) \otimes_{\pi_{\kappa}} \Lambda^{\bullet}(\g)$ and it is always quasi-isomorphic to the module of K\"ahler differential forms $\Omega_{\Lie} (\g)$. Therefore,  the twisting morphism $\varphi$ is a Koszul morphism when $U(\g) \otimes_{\varphi} C \qiso \Omega_{\Lie} (\g)$ is a quasi-isomorphism, or equivalently, by the algebraic fundamental twisting morphism, when $\Lie(C) \qiso \g$ is a quasi-isomorphim or when $C \qiso \Lambda^{\bullet}(\g)$ is a quasi-isomorphism.

\subsubsection{\bf Quadratic setting}

When $\g$ is a quadratic Lie algebra, its Koszul dual coalgebra $\g^{\ash}$ is a co-commutative coalgebra and the Koszul complex is $U(\g) \otimes _{\varkappa} \g^{\ash}$, or equivalently, by the Koszul criterion, when $\Lie(\g^{\ash}) \qiso \g$ is a quasi-isomorphim or when $\g^{\ash} \qiso \Lambda^{\bullet}(\g)$ is a quasi-isomorphism.

\smallskip

This provides examples of quadratic Quillen models for Lie algebras \cite{Quillen3}.

\subsection{Koszul duality theory of quadratic modules over an associative algebra}

\subsubsection{\bf Twisted tensor product for modules}
Let $A$ be an associative algebra. Let $\Po = A$ concentrated in arity $1$, that is $A$ is an associative algebra. A $\Po$-algebra $M$ is a left $A$-module $(M,\, \gamma_{M})$. Dually, $\C = C$ and a $\C$-coalgebra is a left $C$-comodule $(N,\, \Delta_{N})$. Let $\alpha : C \rightarrow A$ be a Koszul morphism. The bar construction on $M$ is the chain complex $C\otimes_{\alpha} M := (C \otimes M,\, d_{\alpha} := d_{C \otimes M} + d_{\alpha}^{r})$ where $d_{\alpha}^{r}$ is the composite of
$$C \otimes M \xrightarrow{\Delta_{C} \otimes id_{M}} C \otimes C \otimes M \xrightarrow{id_{C} \otimes \alpha \otimes id_{M}} C \otimes A \otimes M \xrightarrow{id_{C} \otimes \gamma_{M}} C \otimes M.$$
The cobar construction on $N$ is given dually by $A \otimes_{\alpha} N$.

\smallskip

The twisted tensor product is given by the cobar construction $A \otimes_{\alpha} N$ and a twisting morphism $\varphi : sN \rightarrow N$ is Koszul when $A \otimes_{\alpha} N \qiso A \otimes_{\alpha} C \otimes_{\alpha} M$ is a quasi-isomorphism. Since $\alpha : C \rightarrow A$ is a Koszul morphism, this is equivalent to $A \otimes_{\alpha} N \qiso M$ is a quasi-isomorphism.

\subsubsection{\bf Quadratic module}

Assume now that $A = A(E,\, R)$ is a quadratic algebra $T(E)/(R)$ and $C = A^{\ash}$ its Koszul dual coalgebra. We assume moreover that $A$ is a Koszul algebra, so $\alpha = \kappa$ is a Koszul morphism. An $A$-module $M$ is a \emph{quadratic $A$-module} if $M = (A \otimes V)/A \cdot S$ where $V$ is a vector space of generators, $S \subseteq A\otimes V$ is a subvector space of relations and $A \cdot S$ is the image of the map $A \otimes S \mono A \otimes A \otimes V \xrightarrow{\gamma_{A} \otimes id_{V}} A \otimes V$. We define dually the notion of quadratic comodule and we get the Koszul dual $A^{\ash}$-comodule $M^{\ash}$.

\smallskip

In this case, the Koszul complex is equal to the cobar construction $A \otimes_{\kappa} M^{\ash}$ and the Koszul criterion collapses to $A \otimes_{\kappa} M^{\ash} \qiso M$ is equivalent to $M^{\ash} \qiso A^{\ash} \otimes_{\kappa} M$. Thus we recover the Koszul duality theory of $A$-modules given in \cite{PolishchukPositselski}.

\subsubsection{\bf Modules over a commutative algebra}

Assume now that $A$ is a commutative algebra. In this case, the Koszul duality theory provides resolutions of modules  needed in algebraic geometry \cite{Eisenbud}. When $A = \K [x_{1},\, \ldots ,\, x_{n}]$ or $A = \K [x_{1},\, \ldots ,\, x_{n}]/(I)$ where $I$ is homogeneous of degree $2$, the Koszul dual comodule $M^{\ash}$ of a quadratic module $M$ provides a good candidate for the \emph{syzygies} of $M$.

\subsection{Other fields of applications}

It is also possible to apply the present Koszul duality theory to many other examples of the literature. For example, this Koszul duality theory for algebras applies to the operad encoding Poisson or Leibniz algebras, with possible applications in differential and Poisson geometry \cite{Kosmann, Fresse2}, to the operad encoding PreLie algebras, with possible applications in algebraic combinatorics and links with renormalisation theory in theoretical physics \cite{ChapotonLivernet, ConnesKreimer}, to the operad encoding hyper-commutative or gravity algebras, linked with the Gromov-Witten invariants \cite{Getzler, Manin}.

\smallskip

There are several ways to prove Koszulity for associative algebra that we plan to extend to algebras over an operad such as Poincaré-Birkhoff-Witt bases \cite{Priddy, Hoffbeck} and Gr\"obner bases \cite{Buchberger, DotsenkoKhoroshkin, BokutChenLi}.

\smallskip

Another direction is also to extend this Koszul duality theory beyond the homogeneous quadratic case: when the algebra has quadratic and linear relations, the Koszul dual coalgebra should have an extra differential \cite{Priddy, Ga-CaToVa}, when the algebra has quadratic, linear and constant relations, the Koszul dual coalgebra should have an extra differential and a curvature \cite{PolishchukPositselski, HirshMilles} and when the algebra has quadratic and higher relations, the Koszul dual coalgebra should be a homotopy coalgebra \cite{MerkulovVallette}.

\appendix

\section{Comparison Lemma}\label{comparisonlemma}

In this appendix, we prove a Comparison Lemma for $\Po$-algebras and $\C$-coalgebras generalizing the associative case \cite{Cartan}. We recall that all the chain complexes are non-negatively graded and we assume that $\Po(0)=0$ and that $\C(0)=0$.

\subsection{Weight graded module}\label{weightsetting}

In the sequel, we consider a \emph{bigraded module} $V$, that is a family of modules $\{V_d^{(n)} \}_{n,\, d \geq 0}$. The lower index is the homological degree and the upper one is the weight grading. A \emph{weight graded dg module}, \emph{wdg module} for short, is a bigraded module endowed with a differential which preserves the weight and lowers the homological degree by $-1$. We say that a wdg algebra or wdg coalgebra $V$ is \emph{connected} when it satisfies $V = V^{(1)} \oplus V^{(2)} \oplus \cdots$. Moreover, the structure maps, as $\gamma_{A}$, preserve the weight.

\begin{thmAppendix}[Comparison Lemma]
Let $\alpha : \C \rightarrow \Po$ be an operadic Koszul morphism between a wdg connected cooperad $\C$ and a wdg connected operad $\Po$. Let $A$ be a wdg connected $\Po$-algebra and $C$, $C'$ be wdg connected $\C$-coalgebras. Let $g : sC \rightarrow sC'$ be a morphism of wdg $\Cs$-coalgebras. Let $\varphi : sC \rightarrow A$ and $\varphi' : sC' \rightarrow A$ be two algebraic twisting morphisms, such that $\varphi' \circ g = \varphi$. The morphism $g$ is a quasi-isomorphism if and only if $id_A \otimes^{\Po} g : A\otimes^{\Po}C \qiso A \otimes^{\Po} C'$ is a quasi-isomorphism.
\end{thmAppendix}

\begin{pf}
Since $\Po$, $\C$, $A$ and $C$ are weight graded, the dg modules $\Po \circ (\Po(A),\, C)$ and $\Po(A,\, C)$ are also weight graded. Moreover, the maps $c_{0}$ and $c_{1}$ (see Section \ref{cotangentcomplex}) in the coequalizer
$$\xymatrix{\Po \circ (\Po (A),\, C) \ar@<0.5ex>[r]^{c_{0}} \ar@<-0.5ex>[r]_{c_{1}} & \Po \circ (A,\, C) \ar@{->>}[r]^{proj} & A\otimes^{\Po} C}$$
preserve the weight gradings, so $proj (\Po(A,\, C)^{(n)}) =: (A\otimes^{\Po} C)^{(n)}$ defines a weight grading on $A\otimes^{\Po} C$. We denote by $M = \oplus_{n \geq 0}M^{(n)}$ the wdg $A$-module $A\otimes^{\Po} C$ and by $M' = \oplus_{n \geq 0} M'^{(n)}$ the wdg $A$-module $A\otimes^{\Po} C'$. We define a filtration $F_p$ on $M^{(n)}$ by the formula
$$F_p (M^{(n)}) := \bigoplus_{m+r \leq p} (A\otimes^{\Po} C^{(r)}_m)^{(n)}.$$
The differential on $M$ is given by $d_{\varphi} = d_{A\otimes^{\Po} C} - d_{\varphi}^{l} = d_{A\otimes^{\Po} \K} \otimes id_{C} + id_{A\otimes^{\Po} \K} \otimes d_{C} - d_{\varphi}^{l}$ (see Section \ref{cotangentcomplex} for a definition of $d_{\varphi}^{l}$). We have
$$\left\{ \begin{array}{lcll}
d_{A\otimes^{\Po} \K} \otimes id_{C} & : & F_{p} \rightarrow F_{p} &\\
id_{A\otimes^{\Po} \K} \otimes d_{C} & : & F_{p} \rightarrow F_{p-1} &\\
d_{\varphi}^{l} & : & F_{p} \rightarrow F_{p-2} & \textrm{since $|\alpha| = -1$ and $|\Po| \geq 0$, and $\alpha^{(0)} = 0$.}
\end{array} \right.$$
Thus $F_p$ is a filtration of the chain complex $M^{(n)}$. We denote by $E_{p,\, q}^{\bullet}$ the associated spectral sequence. We have
$$E_{p,\, q}^0(M^{(n)}) = F_p (M^{(n)})_{p+q}/F_{p-1} (M^{(n)})_{p+q} = \bigoplus_{r=0}^{n} (A\otimes^{\Po} C^{(r)}_{p-r})_{p+q}^{(n)} = \bigoplus_{r=0}^{n} (A\otimes^{\Po} \K)_{q+r}^{(n - r)} \otimes C^{(r)}_{p-r}.$$
The study of the differential on $M$ shows that $d^{0} = d_{A \otimes^{\Po} \K} \otimes id_{C}$ and $d^{1} = id_{A\otimes^{\Po} \K}\otimes d_{C}$. Hence
$$E_{p,\, q}^2 (M^{(n)}) =\bigoplus_{r = 0}^{n} \mathrm{H}_{q+r}((A\otimes^{\Po} \K)_{\bullet}^{(n - r)}) \otimes \mathrm{H}_{p-r}(C_{\bullet}^{(r)}).$$
The filtration $F_p$ being exhaustive and bounded below, we can apply the classical theorem of convergence of spectral sequences (Theorem $5.5.1$ of \cite{Weibel}) to get $E_{p,\, q}^{\infty}(M^{(n)}) \cong gr_{p}$H$_{p+q}(M^{(n)})$.

We can define the same filtration on $M'$ and we obtain the same result of convergence of spectral sequences.
\begin{itemize}
 \item When $g$ is a quasi-isomorphism, we get that $E_{p,\, q}^2(M^{(n)}) \xrightarrow{\textrm{H}_{\bullet} (id_{A\otimes^{\Po}\K})\otimes \textrm{H}_{\bullet}(g)} E_{p,\, q}^2({M'}^{(n)})$ is an isomorphism. Since $\varphi' \circ g = \varphi$, the map $\textrm{H}_{\bullet} (id_{A\otimes^{\Po}\K})\otimes \textrm{H}_{\bullet}(g)$ is an isomorphism of chain complexes and the pages $E_{p,\, q}^r(M^{(n)})$ and $E_{p,\, q}^r(M^{(n)})$ are isomorphic for all $r\geq 2$. By the convergence's theorem of spectral sequences, we get that $gr_{p}$H$_{p+q}(M^{(n)}) \cong E_{p,\, q}^{\infty}(M^{(n)}) \xrightarrow{\textrm{H}_{\bullet} (id_A\otimes^{\Po} g)} E_{p,\, q}^{\infty}({M'}^{(n)}) \cong gr_{p}$H$_{p+q}({M'}^{n})$ is also an isomorphism.
\item Assume now that $id_A \otimes^{\Po}g$ is a quasi-isomorphism of dg $A$-modules. Let us work by induction on the weight $n$. When $n=0$, the map $g^{(0)} : 0 \rightarrow 0$ is a quasi-isomorphism. Suppose now that $g^{(n-1)}$ is a quasi-isomorphism. We consider the mapping cone of $f^{(n)} := (id_A \otimes^{\Po}g)^{(n)} : M^{(n)} \rightarrow {M'}^{(n)}$ defined by $cone(f^{(n)}) := s^{-1}M^{(n)} \oplus M'^{(n)}$ and the associated filtration $F_{p}(cone(f^{(n)})) := s^{-1}F_{p}(M^{(n)}) \oplus F_{p}(M'^{(n)})$, which satisfies $E_{\bullet,\, q}^{1}(cone(f^{(n)})) = cone(E_{\bullet,\, q}^{1}(f^{(n)}))$. The mapping cone of $E_{\bullet,\, q}^{1}(f^{(n)})$ fits into a short exact sequence, which induces the long exact sequence
$$\hspace{0.7cm} \cdots \rightarrow \mathrm{H}_{p+1}\left(cone(E_{\bullet,\, q}^{1}(f^{(n)}))\right) \rightarrow \mathrm{H}_{p}\left(E_{\bullet,\, q}^{1}(M^{(n)})\right) \xrightarrow{\mathrm{H}_{p}\left(E_{\bullet,\, q}^{1}(f^{(n)})\right)} \mathrm{H}_{p}\left(E_{\bullet,\, q}^{1}(M'^{(n)})\right)$$
$$\hspace{9.4cm} \rightarrow \mathrm{H}_{p} \left(cone(E_{\bullet,\, q}^{1}(f^{(n)}))\right) \rightarrow \cdots.$$
This induces the long exact sequence $(\xi_{q})$
$$\begin{array}{l} \hspace{0.1cm}(\xi_{q})\
\cdots \rightarrow E_{p+1,\, q}^{2}(cone(f^{(n)})) \rightarrow E_{p,\, q}^{2}(M^{(n)}) \xrightarrow{E_{p,\, q}^{2}(f^{(n)})} E_{p,\, q}^{2}(M'^{(n)}) \rightarrow E_{p,\, q}^{2}(cone(f^{(n)})) \rightarrow \cdots
\end{array}$$
where $E_{p,\, q}^{2}(f^{(n)})$ is given by $\mathrm{H}_{\bullet}(id_{A\otimes^{\Po} \K}) \otimes \mathrm{H}_{\bullet} (g)$. Since $\mathrm{H}_{\bullet}(id_{A\otimes^{\Po} \K})$ is an isomorphism (it is the identity) and $(A\otimes^{\Po} \K)^{(0)} = \K$, the formula for $E_{p,\, q}^2 (M^{(n)})$ given above and the induction hypothesis tells us that $\mathrm{H}_{\bullet}(id_{A\otimes^{\Po} \K}) \otimes \mathrm{H}_{\bullet} (g)$ is an isomorphism, except for $q = -n$ when H$_{q+n}((A\otimes^{\Po} \K)^{(0)}) = \K \neq 0$. The long exact sequence $(\xi_{q})$ for $q \neq -n$ and for all $p$ gives that $E_{p,\, q}^{2}(cone(f^{(n)})) = 0$. Thus, the spectral sequence collapses at rank $2$ and $E_{p,\, q}^{2}(cone(f^{(n)})) = E_{p,\, q}^{\infty}(cone(f^{(n)}))$. Moreover the spectral sequence $E_{p,\, q}^{\bullet}(cone(f^{(n)}))$ converges to $\mathrm{H}_{p+q}(cone(f^{(n)}) = 0$ (since $f^{(n)}$ is a quasi-isomorphism), so $E_{p,\, -n}^{2}(cone(f^{(n)})) = 0$ for all $p$. Finally the spectral sequence $(\xi_{n})$ gives the isomorphism
$$\mathrm{H}_{p-n}(C_{\bullet}^{(n)}) = E_{p,\, -n}^{2}(M^{(n)}) \xrightarrow{\mathrm{H}_{\bullet}(g^{(n)})} E_{p,\, -n}^{2}(M'^{(n)}) = \mathrm{H}_{p-n}({C'}_{\bullet}^{(n)}),$$
for every $p$. Hence, $g^{(n)}$ is a quasi-isomorphism. This prove the result by induction.
\end{itemize}
$\cqfd$
\end{pf}

\section*{Acknowledgments}

I am grateful to Bruno Vallette for his fruitful ideas and his constant help. I thank Éric Hoffbeck and David Chataur for their comments on the preliminary version of this paper.

\bibliographystyle{alpha}
\bibliography{bib}
\ \\

\noindent
\textsc{Joan Millès, Laboratoire J. A. Dieudonné, Université de Nice - Sophia Antipolis, Parc Valrose, 06108 Nice Cedex 02, France}\\
E-mail : \texttt{joan.milles@math.unice.fr}\\
\textsc{URL :} \texttt{http://math.unice.fr/$\sim$jmilles}

\end{document}